\documentclass{article}
\usepackage{times}
\usepackage[letterpaper,top=2cm,bottom=2cm,left=3cm,right=3cm,marginparwidth=1.75cm]{geometry}
\usepackage{hyperref}
\usepackage{amsmath, amsfonts,amssymb,amsthm}
\usepackage{mathrsfs} 
\usepackage{lipsum}
\usepackage{amsfonts}
\usepackage{graphicx}
\usepackage{epstopdf}
\usepackage{algorithmic}
\usepackage{float}
\usepackage{xcolor}
\usepackage{tikz}
\usepackage{cleveref}
\usepackage{float}

\theoremstyle{remark}
\newtheorem{remark}{Remark}
\theoremstyle{definition}
\newtheorem{defn}{Definition}

\usepackage[]{graphicx}

\usetikzlibrary{arrows.meta}
\usepackage{mathtools}

\newcommand{\bs}[1]{\boldsymbol{#1}}

\newcommand{\bOmega}{\bar{\Omega}}
\newcommand{\bGamma}{\bar{\Gamma}}
\newcommand{\bbq}{\bar{\bs{q}}}
\newcommand{\bbx}{\bar{\bs{x}}}
\newcommand{\flowmap}{\bar{\bs{\chi}}}
\newcommand{\bbF}{\bar{\bF}}

\newcommand{\bJ}{\bar{J}}
\newcommand{\bbC}{\bar{\bs{C}}}
\newcommand{\bbu}{\bar{\bs{u}}}
\newcommand{\bbv}{\bar{\bs{v}}}
\newcommand{\bbw}{\bar{\bs{w}}}
\newcommand{\bnabla}{\bar{\nabla}}
\newcommand{\bc}{\bar{\bs{c}}}
\newcommand{\bmom}{\bar{\bs{P}}}
\newcommand{\brho}{\bar{\rho}}
\newcommand{\bB}{\bs{B}}

\newcommand{\bbeta}{\bar{\bs{\eta}}}

\newcommand{\ale}{\bs{\alpha}}

\newcommand{\bF}{{\bs{F}}}
\newcommand{\bu}{{\bs{u}}}
\newcommand{\buale}{\bu}
\newcommand{\bv}{{\bs{v}}}
\newcommand{\bw}{\bs{w}}

\newcommand{\R}{\mathbb{R}}
\newcommand{\calQ}{\mathcal{Q}}
\newcommand{\bcalQ}{\bar{\calQ}}
\newcommand{\calH}{\mathscr{H}}
\newcommand{\bcalH}{\bar{\calH}}
\newcommand{\calR}{\mathscr{R}}
\newcommand{\bcalR}{\bar{\calR}}
\newcommand{\calV}{\mathcal{V}}
\newcommand{\bcalV}{\bar{\calV}}
\newcommand{\calW}{\mathcal{W}}
\newcommand{\bcalW}{\bar{\calW}}

\newcommand{\bbG}{\mathbb{G}}
\newcommand{\bbK}{\mathbb{K}}

\newcommand{\bbM}{\mathbb{M}}
\newcommand{\bbS}{\mathbb{S}}

\newcommand{\trafo}{\mathrm{T}}
\newcommand{\dtrafo}{\mathbb{L}}

\newcommand{\D}{\mathrm{D}}
\newcommand{\der}{\mathbb{D}}

\newcommand{\barx}{\bar{\bs{x}}}
\newcommand{\bx}{\bs{x}}
\newcommand{\dbx}{{\rm d}\bs{x}}
\newcommand{\dbbx}{{\rm d}\bar{\bs{x}}}
\newcommand{\bq}{\bs{q}}
\newcommand{\bC}{\mathbb{C}}
\newcommand{\jale}{J}
\newcommand{\bFale}{\bF}
\newcommand{\bFalei}{\bF_{\irm}}
\newcommand{\Ji}{J_{\irm}}
\newcommand{\bvi}{\bs{v}_i}
\newcommand{\Jf}{J_{\frm}}
\newcommand{\Js}{J_{\srm}}
\newcommand{\rhof}{\rho_{\frm}}
\newcommand{\irm}{{i}}
\newcommand{\srm}{{s}}
\newcommand{\frm}{{f}}
\newcommand{\brm}{{b}}
\newcommand{\krm}{{k}}

\newcommand{\pre}{\pi}
\newcommand{\con}{\bs{c}}
\newcommand{\idty}{\bs{I}_d}

\crefname{section}{section}{sections}

\definecolor{andreacol}{rgb}{0.93, 0.53, 0.18}

\title{Discretisation of Eulerian nonlinear elasticity and diffusion \\ using gradient flows}

\author{Andrea Zafferi \thanks{Freie Universit\"at Berlin, Department of Applied and Computer Science,  Arnimalle 9, 14195 Berlin} \thanks{Corresponding Author E-mail: andrea.zafferi@fu-berlin.de} 
    \and Dirk Peschka \thanks{Weierstrass Institute, Anton-Wilhelm-Amo-Str. 39, 10117 Berlin} \thanks{E-mail: dirk.peschka@wias-berlin.de}}

\begin{document}
\date{}
\maketitle

\begin{abstract}
In this study, we introduce a general energy-based modelling approach for viscous poroelastic materials that feature diffusive transport in both Lagrangian and Eulerian frames. Our research produces refined weak formulations by using the reference map concept within the Eulerian configuration. We propose and implement a novel structure-preserving discretisation strategy, utilising mixed finite element methods. This paper highlights the spatial and temporal numerical convergence of our methods through a comparative analysis of Lagrangian and Eulerian schemes, thereby proving the robustness and usability of our approach. Furthermore, in the context of Eulerian multiphase flow, specifically of the quasi-static Euler-Euler type, our study demonstrates the existence of solitary fluid waves within poroviscoelastic media. 
This energy-based approach forms a basis for a deeper understanding of thermodynamical modelling and corresponding discretisation schemes for coupled poroelasticity, flow, and diffusion.
\end{abstract}




\section{Introduction}
Many physical, biological, and geoscientific phenomena involving viscoelastic materials can be described by thermodynamic models formulated in either Lagrangian or Eulerian coordinates \cite{roubivcek2005nonlinear}. Such models naturally consider diffusion-reaction processes, interfacial effects such as capillarity and phase separation, and external fields such as gravity and electrostatic fields, leading to complex interactions across multiple spatial and temporal scales \cite{gurtin1996generalized,tanaka2000viscoelastic,hong2011modeling,mielke2025model}. When large elastic deformations are present, the resulting partial differential equations are highly nonlinear, and even their rigorous formulation in the spatial (Eulerian) frame poses fundamental mathematical and modelling challenges, \emph{e.g.}, the consistent coupling of elastic and inelastic effects \cite{mielke2025general}, the enforcement of Onsager symmetries \cite{doi2011onsager}, the incorporation of order parameters for phase-field models {\cite{garcke2022viscoelastic,agosti2017cahn,colli2023cahn,lazzaroni2018rate,tornquist2025}}, the treatment of fluid–structure interaction \cite{benevsova2023variational,coutand2005motion}, and the choice of appropriate state variables, such as strain measures, mappings, or stresses \cite{marsden1994mathematical,gurtin2010mechanics}. 
A key feature of most of these nonlinearly coupled models is their energy-variational formulation stated in either a Lagrangian or in an Eulerian frame.

Traditionally, Lagrangian approaches have been used for the mathematical treatment of large deformations, \emph{e.g.},  \cite{roubicek2018thermodynamics,van2024finite,schmeller2023gradient}. Alternatively, Eulerian descriptions are often based on formulations that use covariant Lie-derivatives of deformation gradients ($\bs{F}$), left or right Cauchy green tensor ($\bs{B}=\bs{F}\bs{F}^T$ or $\bs{C}=\bs{F}^T\bs{F}$) or related stress tensors, \emph{e.g.}, \cite{edwards1991non,marsden1994mathematical,mielke2025eulerian,snoeijer2020relationship}. These approaches have  been effectively employed in modelling complex fluid-structure interactions \cite{liu2001eulerian} and extended to nonlinear elastic biological tissues, capturing their inherent compressibility and growth dynamics \cite{wei2023eulerian}. In a similar fashion, regarding geological applications an Eulerian large strain model for porous materials was developed in \cite{roubivcek2023viscoelastodynamics}, where the energetics of the system is obtained as part of the mathematical analysis.
Alternatively, the use of the \emph{reference map} has emerged as a valuable tool for handling finite-strain elasticity within Eulerian frameworks with heterogeneous materials, combining the spatial and material viewpoints  \cite{kamrin2012reference}. Therefore, such an approach might be particularly suitable for geophysical applications and porosity waves. 

Porosity waves are relevant for various natural processes by governing the fluid transport in deformable porous media. Porosity waves typically emerge and evolve due to gradients in the pore pressure, mechanical deformation of the porous solid matrix, and gravity acting on pore fluid and solid. Recent works  highlight their relevance for geophysical \cite{alkhimenkov2024shear,tian2014impact,skarbek2016dehydration} and for biological problems \cite{chiarelli2014poroelastic}. 
In particular \cite{alkhimenkov2024shear} emphasises  the interplay of compaction‐driven fluid flow and plastic effects to generate localised fluid flow 
patterns such as solitary porosity waves.
The classical works by Biot provided the foundational theory of poroelasticity, describing wave propagation through fluid-saturated porous media \cite{biot1956a,biot1956b}, and were later extended by McKenzie \cite{mckenzie1984generation}. 
However, accurately modelling porosity waves poses significant challenges, particularly due to the complexity inherent in coupling fluid-structure interactions and multiphase transport phenomena. 
Small-strain approaches are employed frequently but fail to capture the transition from small deformations of elastic solids to large deformations of soft solids or liquids.

Modelling approaches that are based on thermodynamic principles rely on the definition of thermodynamic potentials such as energy or entropy and additional geometric structures or kinematic constructions involving conservation laws to deduce the evolution of the system. Within these frameworks, we highlight the Hamiltonian \cite{morrison1998hamiltonian}, damped Hamiltonian, and gradient system formulations \cite{otto2001geometry}, which naturally incorporate thermodynamic consistency and systematic treatment of dissipative processes. 
Among such approaches, one can find the GENERIC (General Equation for Non-Equilibrium Reversible-Irreversible Coupling) formalism, introduced by Grmela and \"Ottinger \cite{grmela1997dynamics,ottinger1997dynamics}, used to couple reversible and irreversible processes in isolated systems and later extended to thermoelastic solids and to damped Hamiltonian structure \cite{mielke2011formulation}. 
The development of corresponding structure-preserving discretisations promises to generate robust discrete schemes that also provide stable numerical methods for complex, nonlinearly coupled problems, \emph{e.g.}, see \cite{betsch2019energy,jungel2020minimizing} in the context of GENERIC. Such energy-based discretisation schemes have been also used in the context of poroelasticity, \emph{e.g.}, in \cite{egger2021structure,barnafi2021mathematical,schmeller2023gradient}, where saddle-point problems are a reoccurring theme.

In the spirit of previous works  \cite{peschka2022variational,zafferi2022generic}, we introduce energy-variational   structures for nonlinear poroelastic materials and present a weak formulation of their gradient systems.
The goal is to highlight the origin of saddle-point structures and their natural discretisation in Eulerian and Lagrangian formulations. Therefore, in \Cref{sec:framework} we introduce the basic notions and theoretical concepts underlying energy-based models across different coordinate frames. Most importantly, we present a transformation strategy that establishes a structure-preserving  correspondence between energies and dissipative structures in the Lagrangian and Eulerian frameworks.
In \Cref{sec:elasticity} we develop the kinematics of nonlinear elasticity and apply the framework to energy minimisers in the Lagrangian reference configuration. Using the previously introduced transformation, we derive the corresponding Eulerian formulation and analyse convergence of solutions between and in the two descriptions.
Later, in \Cref{sec:extension_concentration} we extend this approach to gradient systems in poroelasticity for a single solid or fluid phase coupled to a diffusing species, and we demonstrate descending poroelastic waves. We compare fluid-like and solid-like regimes, where in the latter diffusion dominates.
Finally, in \Cref{sec:poro} we formulate a gradient system for hyperelastic porous materials coupled to fluid flow of Euler-Euler type in Eulerian coordinates and emphasise the choice of state variables, energy functionals, and dissipation potentials that define the gradient system based on everything learned before.


\section{Formal structure-preserving variational framework}
\label{sec:framework}
The goal of this section is to introduce the variational setting for an extended gradient flow structure that has practical benefits for numerical discretisation. Then, we present a systematic approach to perform non-trivial change of coordinates (or variables), \emph{e.g.}, from Lagrangian to Eulerian coordinates, and maintain the energetic-variational structure that for a discretisation can be beneficial to maintain energy-robustness of the discrete scheme.
\subsection{Gradient flows}
\label{sec:GS-framework}
We define a gradient structure as a triple $(\calQ, \calH, \calR)$ consisting of a state space $\calQ$, an energy functional $\calH$, and a dissipation $\calR$. The elements of the state space will be denoted by $\bq\in\calQ$ and the associated velocity space is $\calV$ with elements $\bv\in\calV$ or $\dot{\bq}=\partial_t\bq\in\calV$ for the partial time derivative of $\bq$. Usually, $\calQ$ is the subset of a Banach or Hilbert space.
The free energy of the system is a functional $\calH:\calQ\to\mathbb{R}$ and dissipative effects are encoded in the convex dissipation potential $\calR(\bq,\cdot):\calV\to\mathbb{R}$. For some $\bs{\xi}\in\calV^*$ and some $\bv\in\calV$, we denote by $\left\langle\cdot, \cdot\right\rangle_{\calV}:\calV^*\times\calV\to\R$ the canonical dual pairing $\left\langle\bs{\xi},\bv\right\rangle_{\calV}=\bs{\xi}(\bv)$. By applying the \emph{Legendre transform} we can define the dual dissipation potential for any element $\bs{\xi}\in\calV^*$ by $\calR^*(\bq,\bs{\xi}):= \sup_{\bv\in\calV}\left(\left\langle\bs{\xi},\bv\right\rangle_{\calV}- \calR(\bq,\bv)\right)$.
Moreover, if the dissipation $\calR(\bq,\cdot)$ is quadratic, then we write
\begin{equation} \label{eqn:diss_pot}
    \calR(\bq,\bv) = \frac{1}{2}\left\langle\bbG(\bq)\bv,\bv\right\rangle_{\calV}\,,
\end{equation}
with the associated positive symmetric operator $\bbG(\bq):\calV\to\calV^*$. In this case, we call the triple $(\calQ,\calH,\bbG)$ a classical gradient system. The resulting evolution equation has the form
\begin{equation}\label{eqn:GS_equation}
    \bbG(\bq)\dot{\bq} = -\D \calH(\bq)\,, \quad \text{in }\,\,\calV^*\,, 
\end{equation}
in the dual of the velocity space, where $\D \calH(\bq)$ is the \emph{Fr\'echet} derivative of $\calH$ defined by
\begin{equation*}
    \left\langle\D \calH(\bq),\bv\right\rangle_{\calV} := \lim_{h\to 0}\frac{\calH(\bq+h\bv) - \calH(\bq)}{h}\,,\text{ for any }\bv\in\calV\,.
\end{equation*}
Alternatively, equation \eqref{eqn:GS_equation} extends to  arbitrary convex dissipation potentials by  
\begin{align}\label{eqn:weak_gs}
    \left\langle \D_{\dot{\bq}}\calR(\bq,\dot{\bq}), \bv \right\rangle = - \left\langle \D_{\bq}\calH(\bq), \bv\right\rangle_{\calV}\,, \qquad \text{ for any }\bv \in \calV\,.
\end{align}
By applying the definition of $\calR^*$ we introduce the Onsager operator $\bbK(\bq):\calV^*\to\calV$ from the relation $\calR^*(\bq,\bs{\xi})=\frac{1}{2}\left\langle\bs{\xi},\bbG(\bq)^{-1}\bs{\xi}\right\rangle=\frac{1}{2}\left\langle\bs{\xi},\bbK(\bq)\bs{\xi}\right\rangle$. Together with \eqref{eqn:GS_equation}, this results in an alternative evolution equation in the velocity space
\begin{equation}\label{eqn:Onsager_equation}
    \dot{\bq} = -\bbK(\bq)\D\calH(\bq)\,, \quad \text{in }\,\,\calV\,.
\end{equation}

\noindent
A further alternative way to represent the evolution equations \eqref{eqn:GS_equation} and \eqref{eqn:Onsager_equation} is to consider the auxiliary space $\calW$ and a linear operator $\bbM^*(\bq):\calW\to\calV^*$ that provides a representation of the generalised forces $\bs{\xi}\in\calV^*$ in terms of $\bs{\eta}\in\calW$ by the defining relation 
\begin{align}
\label{eqn:defeta}
\bbM^*(\bq)\bs{\eta} = \D\calH(\bq)\,.
\end{align}
Assuming such a function $\bs{\eta}$ exists and exploiting this relation and the  operator $\bbM(\bq):\calV\to\calW^*$, defined such that $\left\langle\bbM^*(\bq) \bw,\bv\right\rangle_{\calV} = \left\langle\bbM(\bq) \bv, \bw\right\rangle_{\calW}$, one can rewrite equation \eqref{eqn:Onsager_equation} as a weak form of a saddle-point problem
\begin{subequations}\label{eqn:MAU}
    \begin{align} 
        k(\bs{\eta},\bw) &+ b(\dot{\bq},\bw) &&= 0 &&&&\text{for all }\bw\in\calW\,,
        \\
        b(\bv,\bs{\eta}) & &&=\left\langle\D\calH(\bq),\bv\right\rangle_{\calV} &&&&\text{for all } \bv\in\calV\,,
    \end{align}
    where $k(\bs{\eta},\bw) = \left\langle \bbM^*(\bq)\bw,\bbK(\bq)\bbM^*(\bq)\bs{\eta}\right\rangle_{\calV}$ is symmetric, positive definite and $b(\bv,\bw)=\left\langle\bbM^*(\bq)\bw,\bv\right\rangle_{\calV}$. 
\end{subequations}
When no confusion is possible, we will also omit the subscript $\calV$, leaving $\left\langle\cdot,\cdot\right\rangle_{\calV}=\left\langle\cdot, \cdot\right\rangle$. Furthermore, we use $(u,v)=\int_\Omega uv\,\mathrm{d}x$ for the scalar, vectorial or tensorial $L^2$ inner product. Otherwise, for a scalar product in a Hilbert space $\calV$ we will use $(\cdot,\cdot)_\calV$.

System \eqref{eqn:MAU} accounts only for dissipative effects coming from the operator $\bbK$ (or $\calR^*)$ represented by the action on $\bs{\eta}\in\calW$. However, this definition can be further extended to account for additional contributions acting on velocities $\bv\in\calV$ by introducing an operator $\bbS(\bq):\calV\to\calV^*$ (in the same spirit as $\bbG$) and an associated symmetric and positive definite bilinear form $s(\bv_1,\bv_2) = \left\langle \bbS(\bq) \bv_1,\bv_2\right\rangle$, 
which is added to the system leading to \begin{subequations}  \label{eqn:MAUS}
    \begin{align}\label{eqn:MAUS-1}
        k(\bs{\eta},\bw) &+ b(\dot{\bq},\bw) &&= 0 &&&&\forall \bw\in\calW\,,
        \\ \label{eqn:MAUS-2}
        b(\bv,\bs{\eta}) &- s(\dot{\bq},\bv) &&=\left\langle\D\calH(\bq),\bv\right\rangle  &&&&\forall \bv \in \calV\,.
    \end{align}
\end{subequations}
Observe that in the case $k(\cdot,\cdot)\equiv 0$, then the equations \eqref{eqn:MAUS} reduce simply to a weak formulation of \eqref{eqn:GS_equation} with $\bbS\equiv\bbG$. 
One can verify that the energy decreases along solutions by testing with the time derivative of the solution $\bv=\dot{\bq}$ and $\bw=\bs{\eta}$ and by subtracting the equations \eqref{eqn:MAUS}
\begin{align}\label{eqn:energy_descent}
        \frac{\rm d}{\mathrm{d}t}\calH\bigl(\bq(t)\bigr)=\left\langle\D\calH(\bq),\dot{\bq}\right\rangle = -(k(\bs{\eta},\bs{\eta}) + s(\dot{\bq},\dot{\bq}))\leq 0 \,.
\end{align}

\begin{remark}[Some Examples]\label{rem:examples}
There are some trivial examples for structures for an evolution law $\dot{\bq}=-\mathbb{K}(\bq)\mathrm{D}\calH(\bq)$. Therefore, for the moment assume that $\calV=\calQ$ is a Hilbert space with scalar product $(\cdot,\cdot)_\calV$ and the operator $\mathbb{K}(\bq):\calV^*\to\calV$ is invertible. Furthermore, 
the Riez isomorphism $\mathbb{R}_\calV:\calV\to\calV^*$ is defined $\langle\mathbb{R}_\calV\bv,\bs{\eta}\rangle=(\bv,\bs{\eta})_{\calV}$.\\

\begin{enumerate}
\item $\calW=\calV$: 
Then $\mathbb{M}=\mathbb{R}_\calV$ and $k(\bs{\eta},\bs{w})=\langle \mathbb{R}_\calV \bs{w},\mathbb{K}(\bq)\mathbb{R}_\calV \bs{\eta}\rangle$ and $s=0$.
\item $\calW=\calV^*$: Then $\mathbb{M}=\mathbb{I}$ and $k(\bs{\eta},\bs{w})=\langle \bs{w},\mathbb{K}(\bq)\bs{\eta}\rangle$ and $s=0$.
\item $\calW=0$: Then $k=0$ and $s(\bv_1,\bv_2)=\langle \mathbb{K}^{-1}(\bq)\bv_1,\bv_2\rangle$.
\item \textbf{Gelfand triple}: Let us assume that $\calW$ is an Hilbert space and that $\calV$ is dense and continuously embedded into $\calW$. Then after identification of $\calW$ with its dual we will have $\calW^*\equiv\calW\subset \calV^*$ also continuously, giving altogether
    $\calV\subset\calW\equiv\calW^*\subset\calV^*$. The triple $(\calV, \calW, \calV^*)$ is called a \emph{Gelfand triple} (sometimes also called \emph{Rigged Hilbert space} or \emph{evolution triple}). This kind of structured spaces has shown to be useful in the analysis and modelling of partial differential equations, since this construction allows us to represent the duality product between $\calV$ and $\calV^*$. 

Then we can define $b:\calV\times\calW\to\mathbb{R}$ via $b(\bv,\bw)=(\bv,\bw)_{\calW}$. However, note that the corresponding operator $\mathbb{M}^*:\calW\to\calV^*$ is not surjective in general. The corresponding bilinear form would be defined $k(\bs{\eta},\bw)=\langle \mathbb{M}^*\bs{\eta},\mathbb{K}(\bq)\mathbb{M}^*\bw\rangle_\calW$.
\end{enumerate}

\end{remark}

\subsection{Structure preserving transformations}
\label{sec:trafo}
One of the main novelties of this work is to present different systematic approaches of mapping the energy-based variational formulations and their mixed finite element discretisations between different coordinate frames, \emph{i.e.}, between Lagrangian to Eulerian frames. Later, this technique will allow us to construct more complex examples relevant for geophysical applications.

Consequently, in this part we discuss how general transformations influence the constitutive components of this formalism. 
First, let us consider two spaces $\bcalQ$ and $\calQ$ and a sufficiently smooth mapping relating the two sets of state variables  $\bq=\trafo(\bbq)$. For a given energy functional $\bar{\calH}:\bcalQ\to\R$ we assume that there exists an energy functional $\calH:\calQ\to\R$ s.t. 
\begin{align}
\label{eqn:closure_energy}
\bar{\calH}(\bbq)=\calH(\trafo(\bbq))\,
\end{align}
for all $\bbq\in\bcalQ$, then we say that the \emph{closure  condition} for energy is satisfied. Consequently, stationary states of the energy can be computed by seeking $\bq\in\calQ$ such that $\mathrm{D}\calH(\bq)=0$, possibly subject to suitable additional constraints.

Next, we want to relate the evolution $\bbq:[0,T]\to\bcalQ$ to the evolution $\bq:[0,T]\to\calQ$. 
Therefore, consider two structures \eqref{eqn:MAUS}, \emph{i.e.}, two state spaces $\bcalQ, \calQ$, two energy functionals $\bcalH,\calH$ defined on the respective state spaces, and consequently six bilinear forms $\bar{b}, \bar{k}, \bar{s}$ and $b, k, s$ each defined on the respective velocity and auxiliary space. 
For the transformation $\trafo:\bar{\calQ}\to\calQ$ we have a linearisation, the so-called Fr\'echet derivative, $\dtrafo(\bbq):=\D_{\bbq}\trafo(\bbq):\bar{\calV}\to\calV$ to map velocities of the original state to velocities of the transformed state, \emph{i.e.}, $\bv=\dtrafo(\bbq)\bbv$. The transformation so introduced can naturally relate the geometric structures involved. There are several alternative ways to derive equivalent formulations on the saddle-point problem. 

Firstly, the dissipation potentials or dual dissipation potentials will also have to satisfy the \emph{closure condition} $\calR(\trafo(\bbq),\dtrafo\bbv)=\bar{\calR}(\bbq,\bbv)$ for all $\bbv\in\bcalV$ or $\calR^*(\trafo(\bbq),\bs{\xi})=\bar{\calR}^*(\bbq,\dtrafo^*\bs{\xi})$ for all $\bs{\xi}\in\calV^*$, respectively. For the map of the bilinear form $s,\bar{s}$ this leaves the only choice to satisfy a \emph{closure condition} for $s$, \emph{i.e.},
\begin{align}
\label{eqn:closure_s}
s(\trafo(\bbq);\dtrafo\bbv_1,\dtrafo\bbv_2)=\bar{s}(\bbq;\bbv_1,\bbv_2)\,,
\end{align}
for all $\bbv_1,\bbv_2\in\bcalV$, at least if we assume that we map $\bar{s}$ to $s$ and do not mix the role of $s$ and $k$, which in principle might be also possible. Note that usually $s,\bar{s}$ depend on the corresponding state $\bq,\bbq$. For the bilinear form $\bar{k},k$ things are slightly more complex, since also the bilinear forms $b,\bar{b}$ are involved. Therefore, let us first introduce the dual adjoint.
\begin{defn}[Dual adjoint with respect to dual systems]
Let \( \bar{X} \) and \( X \) be vector spaces over \( \mathbb{R} \), and let \( L : \bar{X} \to X \) be a linear map.
Assume we are given two dual systems \( (\bar{X}, \bar{Y}, \bar{b}) \) and \( (X, Y, b) \), where
\[
\bar{b} : \bar{X} \times \bar{Y} \to \mathbb{R}, \quad b : X \times Y \to \mathbb{R}
\]
are the bilinear maps that are associated with the pairings.

Then the \emph{dual adjoint of \( L:\bar{X}\to X \) with respect to the dual systems}, denoted by
\[
L^{*,b} : Y \to \bar{Y},
\]
is defined as the (unique) linear map satisfying
\[
b(L \bar{x}, y) = \bar{b}(\bar{x}, L^{*,b} y)
\quad \text{for all } \bar{x} \in \bar{X}, \; y \in Y.
\]
\end{defn}

\paragraph*{Method 1}
With this definition, let us assume the spaces $\calV,\bcalV,\calW,\bcalW$ and the corresponding bilinear forms (dual systems) $b,\bar{b}$ are already given. Then we say that $k,\bar{k}$ satisfy the \emph{first closure condition} for $k$ if
\begin{align}
\label{eqn:closure_k1}
k(\trafo(\bbq);\bs{\eta}_1,\bs{\eta}_2)=\bar{k}(\bbq;\bbeta_1,\bbeta_2)\,,
\end{align}
where $\bbeta_i=\dtrafo^{*,b}\bs{\eta}_i$ for all $\bs{\eta}_i\in\calW$ for $i=1,2$ using $\dtrafo^{*,b}$ the dual adjoint with respect to the dual systems. This is equivalent to the closure condition for the dual dissipation potential in the sense that, using the chain rule $\mathrm{D}\bcalH(\bbq)=\dtrafo^*(\bbq)\mathrm{D}\calH(\bq)$, we obtain 
\begin{align*}
\langle\mathrm{D}\bcalH(\bbq),\bbv\rangle = \bar{b}(\bbv,\bar{\bs{\eta}}) = \langle \dtrafo^*(\bbq)\mathrm{D}\calH(\bq),\bbv\rangle=\langle \mathrm{D}\calH(\bq),\dtrafo(\bbq)\bbv\rangle = b(\dtrafo(\bbq)\bbv,\bs{\eta})\,,
\end{align*}
which is consistent with $\bbeta_i=\dtrafo^{*,b}\bs{\eta}_i$ above implies \eqref{eqn:closure_k1} and thus
\begin{align*}
\frac{\rm d}{\mathrm{d}t}\calH(\bq)=\langle\mathrm{D}\calH(\bq),\dot{\bq}\rangle = -\bigl(s(\dot{\bq},\dot{\bq})+k(\bs{\eta},\bs{\eta})\bigr)=-\bigl(\bar{s}(\dot{\bbq},\dot{\bbq})+\bar{k}(\bar{\bs{\eta}},\bar{\bs{\eta}})\bigr)\,,
\end{align*}
for any solution. This methodology has been used in \cite{peschka2022variational} for damped Hamiltonian systems, where $b,\bar{b}$ are the $L^2$ inner products in current and referential coordinates. This allowed us to map canonical structures that appeared in mechanical problems to the Eulerian structures, obtaining Lie-Poisson structures that already have been derived by Morrison \cite{morrison1998hamiltonian} with a slightly different methodology. 

\paragraph*{Method 2} Suppose only $\bar{b}$ and $\trafo,\dtrafo$ are given but $b$ is unknown and to be determined. Furthermore, assume that a space and operator $\mathbb{L}_\calW:\bcalW\to\calW$ is given. Then we say $k,\bar{k}$ satisfy the \emph{second closure condition} for $k$ if
\begin{subequations}
\label{eqn:closure_k2}
\begin{align}
k(\trafo(\bbq);\bs{\eta}_1,\bs{\eta}_2)=\bar{k}(\bbq;\bbeta_1,\bbeta_2)\,,
\end{align}
where $\bs{\eta}_i=\mathbb{L}_\calW\bar{\bs{\eta}}_i$ for $i=1,2$ and all $\bar{\bs{\eta}}_i\in\bcalW$. However, to obtain compatibility with the chain rule, we simultaneously need to satisfy the \emph{second closure condition} for $b$
\begin{align}
\bar{b}(\bbv,\bar{\bs{\eta}}) = b(\dtrafo(\bbq)\bbv,\mathbb{L}_\calW\bar{\bs{\eta}})\,,
\end{align}
\end{subequations}
which jointly imply the equivalence of the dissipation potentials.

In summary, in certain cases method 1 or method 2 might be preferable over the other to obtain a closure condition for mapping between energy variational formulations that have the structure of \eqref{eqn:MAUS}. Later, we will observe that Eulerian formulations of gradient flows with poroviscoelasticity and convective derivatives favour method 2 to incorporate certain convective derivatives into the mixed bilinear form $b$. This could potentially lead to a specific class of mixed (finite-element) weak formulations for Eulerian problems with convective derivatives. In the rest of the paper we are going to present mapped problems coupling nonlinear elasticity and diffusive transport from Lagrangian to Eulerian frame and systematically derive mixed finite element formulations that make use of this saddle point structure.

\section{Dissipative hyperelastic solid}
\label{sec:elasticity}
In this section we discuss some generic Lagrangian and Eulerian formulations for nonlinear elasticity to point out the equivalence of both formulations and to show how the closure conditions introduced before are going to be satisfied. After the introduction of the kinematics, we consider a static problem and a gradient flow that is  complemented with a $L^2$-type dissipation potential. Firstly the variational structure is phrased on the Lagrangian (reference) domain $\bOmega$ and subsequently mapped to the Eulerian (current) domain $\Omega$ with the help of the reference map.

\subsection{Kinematics of Lagrangian and Eulerian variables}
\label{sec:kinematics}
We are concerned with functions defined on domains that deform or with particles at an initial position $\barx$ at time $t=0$ that move in space and are at the location $\bs{x}$ at time $t>0$. 
Throughout this paper we will assume, for simplicity, that the shape of the deformed domain $\Omega(t)$, which is obtained from the initial domain $\bOmega$ using the \emph{flow map} $\flowmap:[0,T]\times\bOmega\to\mathbb{R}^d$ via $\Omega(t)=\flowmap(t,\bOmega)$ never changes, so that we can set $\Omega(t)=\bOmega\equiv \Omega$ from hereon. The trajectory of a material point in this description is $\bs{x}(t)=\flowmap(t,\barx)$. In general we call a material-focused description using variables and functions defined in ${\bOmega}$ a Lagrangian formulation and a position-focused description using variables and functions defined in $\Omega(t)$ an Eulerian formulation. We start by introducing a concise notion for functions (densities) and maps (deformation) in Lagrangian and Eulerian variables.

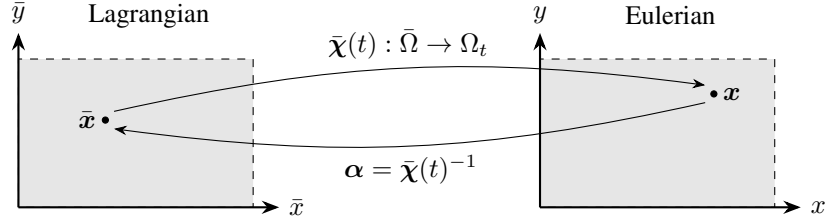
\begin{figure}[H] \label{fig:flow-reference-map}
\centering
\begin{tikzpicture}[scale=1.15,>=Stealth]
    \draw[dashed, fill=gray!20] (-1,0) rectangle (1.7,1.7);
    \draw[thick,->] (-1,0) -- (2,0) node[right] {$\bar{x}$};
    \draw[thick,->] (-1,0) -- (-1,2) node[above] {$\bar{y}$};
    \node at (0.5,2.2) {Lagrangian};
    
    \draw[dashed, fill=gray!20] (5,0) rectangle (7.7,1.7);
    \draw[thick,->] (5,0) -- (8,0) node[right] {$x$};
    \draw[thick,->] (5,0) -- (5,2) node[above] {$y$};
    \node at (6.5,2.2) {Eulerian};

    \draw[->,bend left=10] (0.1,1.1) to node[above] {$\flowmap(t):\bOmega\to\Omega_t$} (6.9,1.4);
    \draw[->,bend left=10] (6.9,1.2) to node[below] {$\bs{\alpha}=\flowmap(t)^{-1}$} (0.1,0.9);

    \filldraw [black] (0,1) circle (1pt) node[anchor=east] {${\bar{\bs{x}}}$};
    \filldraw [black] (7,1.3) circle (1pt) node[anchor=west] {$\bs{x}$};
\end{tikzpicture}
\caption{Flow map $\flowmap$ from Lagrangian to Eulerian configuration and its reference map $\bs{\alpha}=\flowmap^{-1}$.}
\end{figure}

The flow map defines a deformation gradient $\bbF:[0,T]\times\bOmega\to\mathbb{R}^{d\times d}$ via
\begin{align} \label{eqn:def-bbf}
\bbF(t,\bar{\bs{x}})=\frac{\partial \bs{x}}{\partial\bar{\bs{x}}}=\bar{\nabla}\flowmap(t,\bar{\bs{x}}),
\end{align}
with initial data $\flowmap(t=0,\bar{\bs{x}})=\bar{\bs{x}}$ and therefore $\bbF(t=0,\bar{\bs{x}})=\idty$. Local volumetric changes in the domain are computed using its Jacobian 
\begin{equation*}
    \bJ(t,\bar{\bs{x}}):=\det \bbF(t,\bar{\bs{x}})\,. 
\end{equation*}
The main variable that we will be using for describing motion of particles in systems is the displacement that relates to the flow map through the equation
\begin{equation} \label{eqn:def-bbu}
    \bbu(t,\bar{\bs{x}}):=\flowmap(t,\bar{\bs{x}}) - \bar{\bs{x}}\,.
\end{equation}
From this definition we get the following relations with partial derivatives of the displacement $\bbF(t, \barx) = \idty + \bnabla\bbu(t, \barx)$ and $\partial_t \bbu(t, \barx) = \partial_t \flowmap(t, \barx)$. Additionally, we introduce the \emph{reference map} as the inverse of the flow map at time $t$, \emph{i.e.}, $\ale(t,\cdot)=\flowmap^{-1}(t,\cdot)$ or
\begin{align*}
\ale(t,\flowmap(t, \barx))=\barx,
\end{align*}
for all $(t, \barx)\in [0,T]\times\bOmega$, see  \Cref{fig:flow-reference-map}. By using $\ale$ we can get an Eulerian expression for the displacement, that we are going to use as the primary unknown in  \Cref{sec:Eul-hyper-el}
\begin{equation} \label{eqn:def-buale}
    \bbu(t,\ale(t,\bx))=\bx - \ale(t,\bx)=:\bu(t,\bx)\,.
\end{equation}
This allows us to make the following key definition that relates time-dependent Lagrangian functions defined on $\bOmega
$ and Eulerian functions defined on $\Omega$. 
\begin{defn}[Corresponding Lagrangian and Eulerian functions]
\label{def:mappedfunctions}
Let $\bar{a}(t,\bbx)$ an arbitrary (scalar, vectorial, tensorial) function defined for times $t\in(0,T)$ and in space $\bbx\in\bOmega$ and $\flowmap:[0,T]\times\bOmega\to\Omega$ an (invertible) flow map. 
Then we define a \emph{corresponding Eulerian function} $a(t,\bx)$ for times $t\in(0,T)$ and space $\bx\in\Omega$ via
\begin{align*}
a(t,\bx) = \bar{a}(t,\flowmap^{-1}(t,\bx)),
\end{align*}
and vice versa, for any function $a(t,\bx)$ the \emph{corresponding Lagrangian function} $\bar{a}(t,\bbx)$ is
\begin{align*}
\bar{a}(t,\bbx) = a(t,\flowmap(t,\bbx)).
\end{align*}
\end{defn}

\begin{remark}
The \Cref{def:mappedfunctions} is compatible with the definition of displacements $\bbu$ and $\bu$ above and the common definition of Lagrangian and Eulerian velocity $\bar{\bs{v}}(t,\bbx)=\partial_t\flowmap(t,\bbx)$ and $\bs{v}(t,\bx)=\bar{\bs{v}}(t,\flowmap^{-1}(t,\bx))$. For the flow map and the return map itself, this definition produces the trivial identities $\bs{\chi}(t,\bx)=\bx$ and $\bar{\bs{\alpha}}(t,\bbx)=\bbx$. Note that for tensors this definition might be incompatible with a covariant transformation behaviour. 
\end{remark}

Upon differentiation of \eqref{eqn:def-buale} and the identity $\ale(t,\flowmap(t, \barx))=\barx$ we obtain
\begin{subequations}
    \begin{align}\label{eqn:bu-time-der}
        \partial_t \bu(t,\bx) =& -\partial_t \ale(t,\bx)\,,
        \\ \label{eqn:ale-time-der}
        \partial_t \ale(t,\bx) =& -(\nabla \ale)\bv(t,\bx)\,,
        \\
        \bF(t,\bx) =& (\nabla \ale)^{-1}(t,\bx)\,,
    \end{align}
    where $\bv$ denotes the Eulerian velocity field $\bv:= \partial_t\flowmap\circ\ale$. By adding \eqref{eqn:bu-time-der} and \eqref{eqn:ale-time-der} we get 
    \begin{align}
      \bF(t,\bx)\partial_t \bu(t,\bx) =& \bv(t,\bx)\,,
    \end{align}
    where the Eulerian deformation gradient $\bF$ can also be represented as
    \begin{align}
    \bF(t,\bx)= (\nabla \ale)^{-1}(t,\bx) = (\idty - \nabla \bu(t,\bx))^{-1}\,,
    \end{align}
Differentiating the identity $\bu(t,\flowmap(t,\bbx))=\bbu(t,\bbx)$ with respect to time we get the relation
    \begin{align} \label{eqn:disp_vel}
        \partial_t\bu = \left(\idty - \nabla \bu\right)\partial_t\bbu\,.
    \end{align}
    We also state some chain rules helpful in the computation of gradients and determinant of $\bF$ and $\jale$, cf. \cite[eq.(62)]{petersen2008matrix}, 
    \begin{align}
        \der_{\bu}\jale =&\jale \bFale^{\top}:\nabla \Box\,, \\
        \der_{\bu}f(\bFale) =&\bFale^{\top}\left(\partial_{\bFale}f(\bFale)\right)\bFale^{\top}:\nabla\Box\,.
    \end{align}
\end{subequations}
Although the primary unknown for the Lagrangian formulation is the flow map $\flowmap$ and for the Eulerian formulation it is the reference map $\bs{\alpha}$, we will use the conventional displacements $\bbu$ and $\bu$ in the presentation and the computations. This results in a more uniform representation of the nonlinear elasticity in both coordinate frames.

\subsection{Energetics of hyperelastic materials}
We consider a two dimensional bounded domain $\bOmega\subset\mathbb{R}^d$ and the state space $\bcalQ$ consisting of displacements $\bbu:\bOmega\to\mathbb{R}^d$, where we focus on the presentation of the case of dimension $d=2$. As a free energy  of the system we use a functional $\bar{\calH}:\bar{\calQ}\to\R$ that has two additive local energy densities: $\bar{H}_{\rm elastic}$ accounts for the stored mechanical energy of the elastic material and $\bar{H}_{\rm external}$ accounts for external potential energies, \emph{e.g.}, gravity.
The functional is of the specific form
\begin{subequations}
\label{eqn:lag-elastic-energy}
    \begin{align}
        \bar{\calH}(\bbu) =& \int_{\bOmega}\bar{H}_{\rm elastic}(\bbF) + \bar{H}_{\rm external}(\bbu, \bbx)\,{\rm d}\bbx,
        \\
        \text{where\,\,\,\,}&\bar{H}_{\rm elastic}(\bbF) = \frac{\mu}{2}\left(\mathrm{tr}(\bbC - \idty) - 2\log\left(\bJ\right)\right) + \frac{\lambda}{2}(\bJ-1)^2,
        \\
        \text{and\,\,\,\,}&\bar{H}_{\rm external}(\bbu, \bbx) = \bar{\bs{f}}(\bbx)\cdot(\bbx+\bbu) \,,
    \end{align}
\end{subequations}
where $\bbC=\bbF^\top \bbF$ is the right Cauchy-Green tensor and $\lambda$ and $\mu$ are the bulk and shear modulus of the material. Materials with this specific energy density are  referred to as \emph{compressible Neo-Hookean} materials, cf. \cite{attard2003}. In order to account for gravitational effect, the function $\bar{\bs{f}}=\bar{\bs{f}}(\bbx)$ could be taken to be, \emph{e.g.}, $\bar{\bs{f}}=\bar{\rho}\bs{g}$, where $\bar{\rho}\in\R$ is the referential mass density of the material and $\bs{g}\in\R^d$ is a constant vector pointing in the $y$-direction.  This functional produces the following driving force
\begin{align}
\label{eqn:elasticforce}\left\langle\D\bcalH(\bbu),\bbv\right\rangle = \int_{\bOmega}\left(\mu\left(\bbC - \idty\right)\bbF^{-\top} + \lambda(\bJ - 1)\bJ\bbF^{-\top}\right):\bnabla\bbv + \bar{\bs{f}}\cdot\bbv\,\dbbx\,.
\end{align}

With these definitions we can introduce an equivalent Eulerian energy by setting 
\begin{subequations}
\label{eqn:eul-elastic-energy}
    \begin{align}
        {\calH}(\bu) =& \int_{\Omega}{H}_{\rm elastic}(\bF) + {H}_{\rm external}(\bu, \bx)\,{\rm d}\bx,
        \\
        \text{where\,\,\,\,}&{H}_{\rm elastic}(\bF) = \frac{1}{J}\, \bar{H}_{\rm elastic}(\bF)
        \text{\,\,\,\,and\,\,\,\,}{H}_{\rm external}(\bu, \bx) = \frac{1}{J}\,\bs{f}(\bx)\cdot\bx \,,
    \end{align}
\end{subequations}
The energy in \eqref{eqn:lag-elastic-energy} and \eqref{eqn:eul-elastic-energy} satisfy the closure condition \eqref{eqn:closure_energy} trivially as the mapping of the state variable $\bbu$ to $\bu$ defined in \eqref{eqn:def-buale} is invertible as long as the flow map $\flowmap$ is invertible and we define $\calH(\bu)=\bcalH(\trafo^{-1}(\bu))$ via the change of variables formula.

\begin{figure}
\centering
\includegraphics[width=0.41\textwidth]{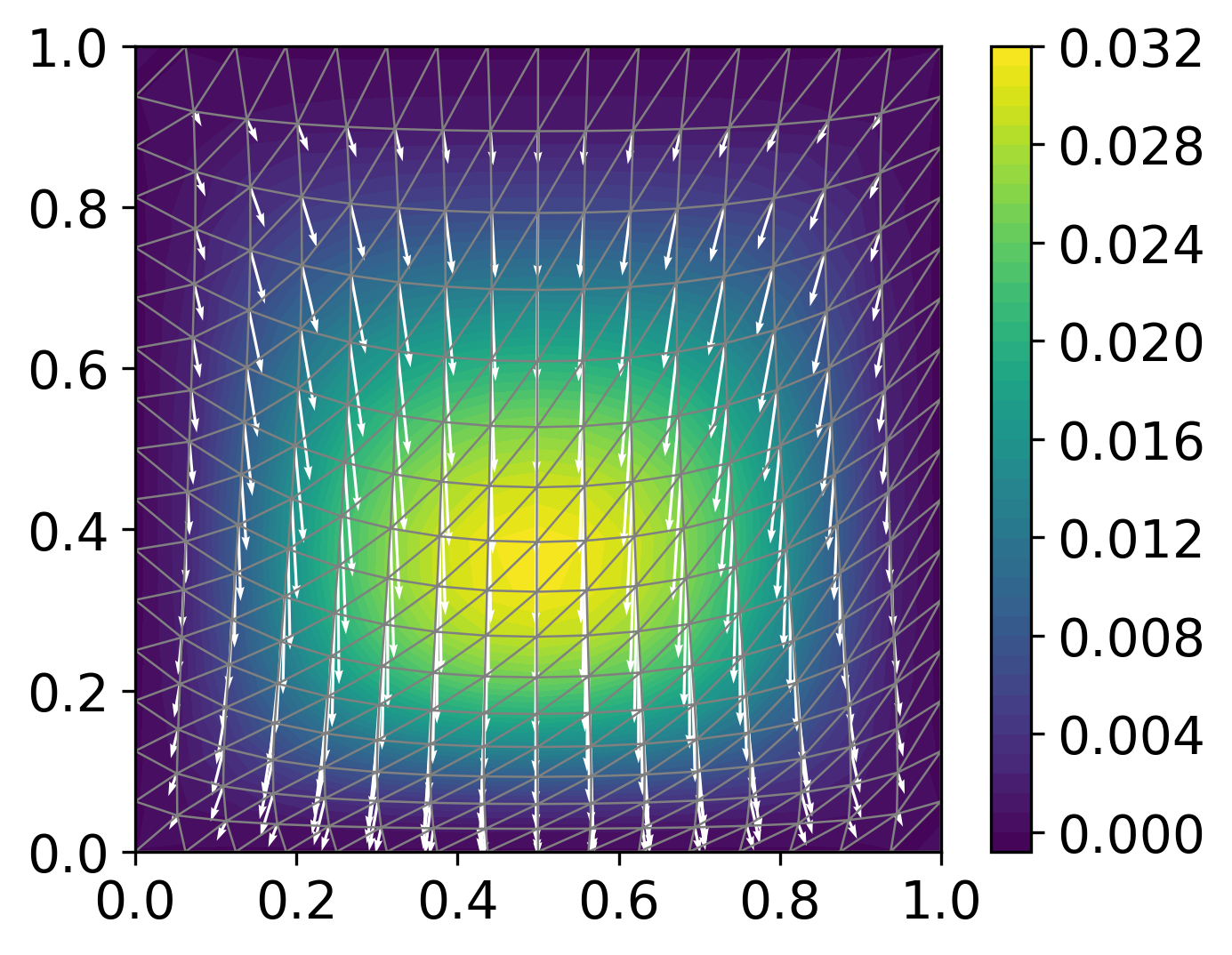}\qquad\qquad
\includegraphics[width=0.41\textwidth]{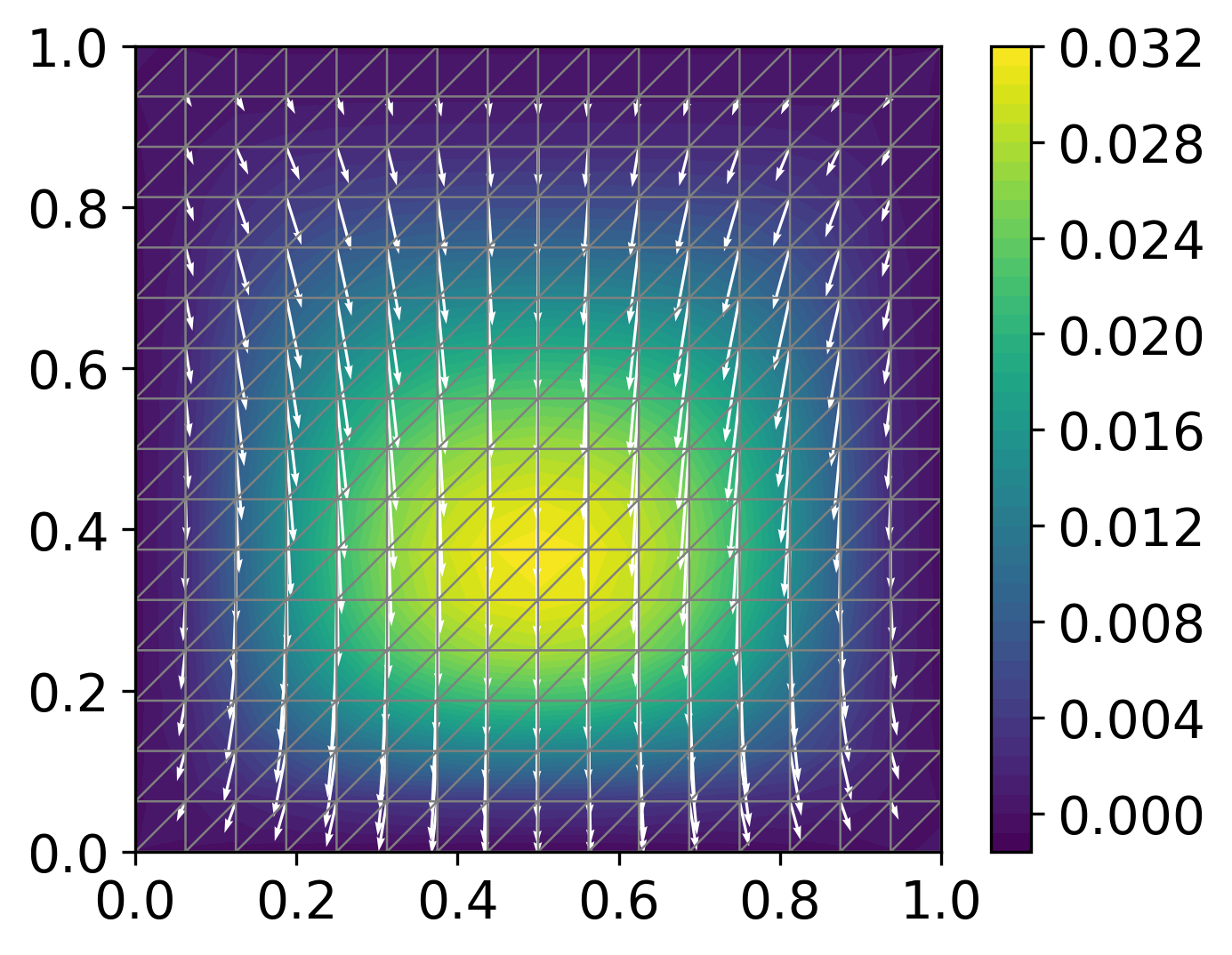}
\caption{(Left) Lagrangian energy minimizer $\bbu$ plotted on a deformed mesh compared to (right) Eulerian energy minimizer $\bu$ on the actual mesh. Shading shows the squared deformation $|\bbu|^2$ and $|\bu|^2$, the white arrows show the vector fields $\bbu$ and $\bu$ and the gray lines show the corresponding $16\times 16$ meshes. The model parameters are $\mu=1$, $\lambda=2$, $\bar{\bs{f}}=(0,6)$}
\label{fig:statsol}
\end{figure}

\begin{figure}
\centering
\includegraphics[width=0.48\textwidth]{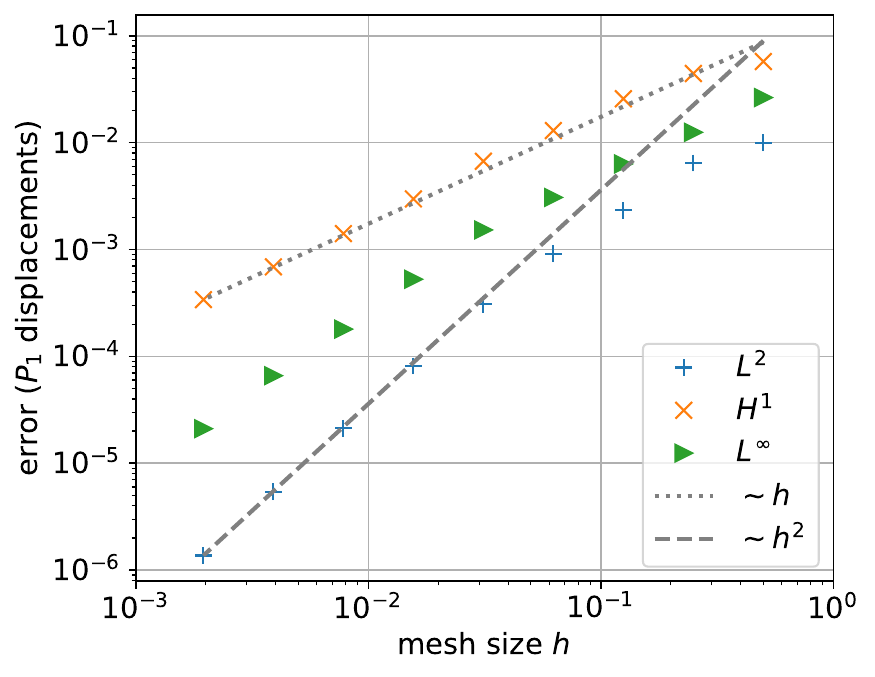}\hfill
\includegraphics[width=0.48\textwidth]{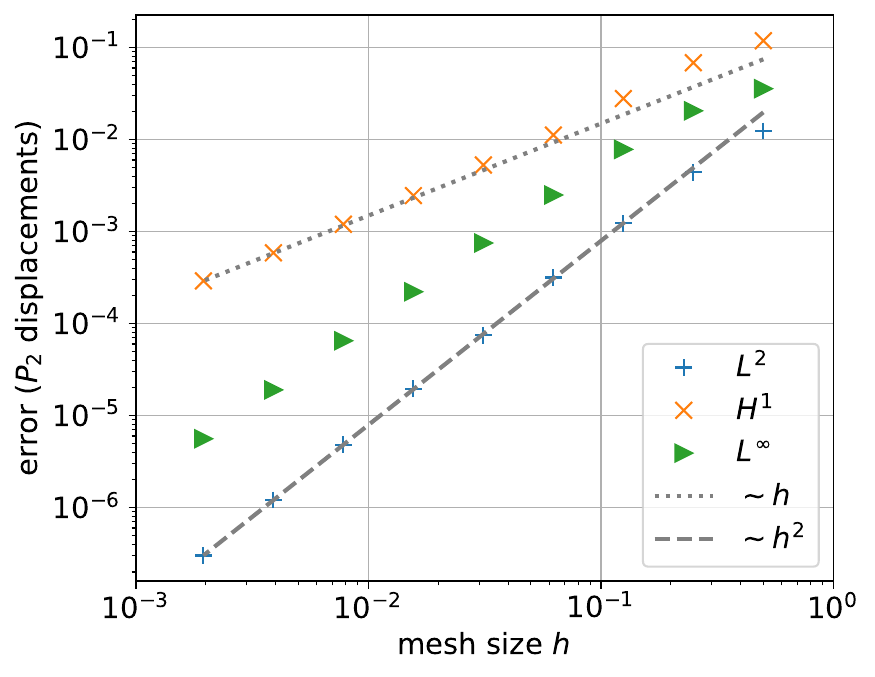}
\caption{Convergence of the error between Lagrangian and Eulerian solution by plotting the norm of the difference $\|\bu_h - \mathrm{proj}_{L^2}\bbu_h\circ\flowmap^{-1}\|$ as a function of mesh size $h=1/(N-1)$ for $N\times N$ meshes for $N=2,...,512$. (Left) shows the convergence for $P_1$ elements and (Right) shows the convergence for $P_2$ elements.}
\label{fig:statsolconv}
\end{figure}

The simplest partial differential equation for such an elastic material is the computation of stationary points (supposed energy minimizers) $\bbu$ and $\bu$ by requiring $\D\bar{\calH}(\bbu)=0$ or $\D\calH(\bu)=0$, respectively.  While both functions describe the same displacement, a direct comparison of solutions is not possible as they are defined on different domains $\bbu:\bOmega\to\R^d$ and $\bu:\Omega\to\R^d$. However, a comparison is possible when projecting on displacement to the other domain, e.g. via $\widehat{\bar{\bu}}(\bx)=\bar{\bu}\circ\flowmap^{-1}$. For discrete solutions $\bar{\bu}_h\in\bar{\calQ}_h$ the corresponding $\widehat{\bar{\bu}}_h$ will usually not lie in the discrete function space of $\bu_h\in\calQ_h$, so that $\widehat{\bar{\bu}}_h$ needs to be projected onto the corresponding solution space via an $L^2$-projection $\mathrm{proj}_{L^2}(\widehat{\bar{\bu}}_h)$. Then we can compute the modelling error between discrete solutions of the Lagrangian and Eulerian formulation via
\begin{align*}
\mathrm{error}(h) = \|\bu_h - \mathrm{proj}_{L^2}(\widehat{\bar{\bu}}_h)\|,
\end{align*}
using a suitable function space norm on $\calQ_h$, e.g. $L^2(\Omega;\R^d)$, $H^1(\Omega,\R^d)$ or $L^\infty(\Omega,\R^d)$.

Representative minimizers for both problems are shown in \Cref{fig:statsol}, where the Lagrangian solution in the left panel is plotted on a deformed mesh to allow for a comparison with the Eulerian solution in the right panel. 
Notice relatively strong deformations with $\max_{\bbx}|\bar{\nabla}\bbu(\bbx)|\approx 1.3$ in the finite strain regime visible in the strong deformation of the computational mesh in the left panel.
In particular the convergence plot shown in \Cref{fig:statsolconv} shows that with both $P_1$ and $P_2$ finite elements, the expected convergence of solution in both cases coincides with the theoretical convergence rate to the exact solution of elliptic problems, \emph{i.e.}, quadratic convergence in the $L^2$ norm and linear convergence in the $H^1$ norm. 
To keep the discussion simple, we do not separately discuss the approximation error due to the projection, which could also contribute to the error for higher-order schemes and sufficiently regular solutions.
This clearly demonstrates that the Lagrangian formulation with flow maps and the Eulerian formulation with reference maps can both be used to effectively solve problems in nonlinear hyperelasticity.

The numerical discretisation shown above is implemented in the finite element framework FEniCS \cite{logg2012automated}, where the functionals are directly provided as in \eqref{eqn:lag-elastic-energy}  and \eqref{eqn:eul-elastic-energy} and then differentiated automatically to solve for stationary solutions.

\subsection{Viscoelasticity with Darcy-type Dissipation}
In the following we introduce nonlinear viscoelasticity in the Lagrangian and the Eulerian frame using the reference map approach. The goal of this example is to introduce the correspondence of these two formulations and show basic (experimental) convergence rates, which should demonstrate the practical usability of the reference map approach derived using gradient flow techniques.
\paragraph{Lagrangian formulation} 
Now we consider the viscous relaxation of a hyperelastic material using a gradient system $(\bcalQ,\bcalH,\bcalR)$.
For the dissipation potential $\bcalR(\bbu,\dot{\bbu})$ we consider a $L^2(\bOmega)$-type dissipation that resembles a Darcy contribution, cf. \cite{zafferiPorousmediaModelReactive2023}. We underline that the viscous stresses arising from these types of potential for single phase systems do not satisfy the Galileian invariance principle, \emph{i.e.}, the value of the stress changes for an observer in constant motion. However it is useful to discuss it here since in Section \ref{sec:poro} we will present a two-phase system and thus consider a Darcy-type flux. The dissipative mechanisms of the system are governed by
\begin{equation}\label{eqn:lag_darcy_dissipation}
\bar{\calR}_\textrm{Darcy}(\bbu,\bbw) := \int_{\bOmega}\frac{\kappa(\bbu)}{2}\vert \bbw\vert^2\,\bJ\,{\rm d}\bbx\,,
\end{equation}
with $\kappa(\cdot)>0$ being the permeability. By taking its Fr\'echet derivative, we define the bilinear form $\srm$
\begin{align}\label{eqn:lag-darcy-diss}
    \srm(\bbw_1, \bbw_2) = \left\langle\bbS(\bbu)\bbw_1, \bbw_2\right\rangle:= \int_{\bOmega} \kappa(\bbu)\bbw_1\cdot \bbw_2 \bJ\,\mathrm{d}\bbx\,\quad \text{ for }\bbw_1,\bbw_2\in\bcalV.
\end{align}
After collecting everything for the abstract structure \eqref{eqn:MAUS} we obtain the following time evolving problem in weak form
    \begin{align*}
         \left(\kappa(\bbu)\dot{\bbu}\bJ,\bbw\right)=-\langle \D\bar{\calH}(\bbu),\bbw\rangle &&\text{for all }\bbw \in \bar{V}\,,
    \end{align*}
where the force $\D\bar\calH$ is  computed in \eqref{eqn:elasticforce}.
The system is completed with homogeneous Dirichlet conditions $\bbu=0$ on the boundary $\bGamma = \partial\bOmega$ and with the initial datum $\bbu_0=0$. 

\paragraph{Eulerian formulation}
\label{sec:Eul-hyper-el}

In the actual configuration, we formulate the system in terms of the Eulerian displacement $\buale(t,\bx)=\bbu(t,\ale(t,\bx))\in\calQ$. Equation \eqref{eqn:def-buale} defines then a mapping $\trafo:\bcalQ\to\calQ$ that relates the Lagrangian displacement $\bbq=\bar{\bs{u}}$ and the Eulerian displacement
\begin{align*}
    \bq = (\buale:\Omega\to\R^d) = \trafo(\bbq)\,.
\end{align*}
Based on the comments of  Section~\ref{sec:trafo}, we map the energy 
\eqref{eqn:lag-elastic-energy} to the current domain $\Omega$ and obtain \eqref{eqn:eul-elastic-energy}. 
This results in the following driving force 
\begin{align*}
    \left\langle\D\calH(\bu),\bw\right\rangle = \int_{\Omega}\bF^{\top}\partial_{\bF}\left(\frac{\bar{H}_{\rm elastic}}{\jale}\right)\bF^{\top}{:}\nabla\bw\,\dbx - \int_{\Omega}\frac{\bs{f}\cdot\bx}{\jale}\bF^{\top}{:}\nabla\bw\,\dbx\,.
\end{align*}
From linearising $\trafo$ we find $\D\mathrm{T}(\bbq)=\dtrafo:\bcalV\to\calV$ using \eqref{eqn:disp_vel} and get $\partial_t\bu = \bF^{-1} \partial_t\bbu$. 
Since $\bF$ is invertible, we can satisfy the closure condition \eqref{eqn:closure_s} for $s$ by setting 
$s(\bv_1,\bv_2)=\bar{s}(\bF\bv_1,\bF\bv_2)$
and using the change of variables formula to obtain the following dissipation potential
\begin{equation} \label{eqn:eul_darcy_diss}
\calR_\textrm{Darcy}(\bu, \bw):=\int_{\Omega}\frac{\kappa(\bu)}{2}\vert \bFale \bw\vert^2{\rm d}\bx\,,
\end{equation}
that produces the bilinear form
\begin{equation*}
    \srm(\bw_1,\bw_2)=\left(\kappa (\bu_{\ale})\bFale\bw_1, \bFale\bw_2\right)\,\quad \text{ for }\bw_1,\bw_2 \in \calV\,.
\end{equation*}
Having specified the properties and the constitutive laws of the material the final evolution equations follow the structure \eqref{eqn:weak_gs}:
    \begin{align}\label{eqn:pde_eul_elastic_darcy}
        \left(\kappa(\bu)\bFale\dot{\bu},\bFale\bw\right) =& -\left\langle\D\calH(\bu),\bw\right\rangle   &&\text{for any }\bw \in  \calV\,.
    \end{align}
    In strong form with boundary and initial conditions the evolution problem reads
    \begin{align*}
        \kappa(\bu)\bF\dot{\bu}\bF^{\top}=&\nabla\cdot\bs{\sigma} - \nabla\cdot\left(\frac{\bs{f}\cdot\bx}{\jale}\bF^{\top}\right)  \, &&\text{in } (0,T]\times\Omega
        \\
        \bu =& \,\bs{0} &&\text{in } (0,T]\times\Gamma\,,
        \\
        \bu =& \,\bs{0}\, &&\text{in } \{0\}\times\Omega\,.
    \end{align*}
    where $\bs{\sigma}:= \jale^{-1}\partial_{\bFale}\bar{H}_{\rm elastic}\bFale^{\top}=\mu\jale^{-1}(\bB - \idty)+ \lambda(\jale -1 )\idty$ is the elastic Cauchy stress tensor, with $\bB=\bF\bF^{\top}$ being the \emph{left Cauchy-Green} deformation tensor or \emph{Green strain tensor}. 
    Equivalently, we can express the elastic Cauchy stress  $\bs{\sigma}=\partial_{\bF}(H_{\rm elastic})\bF^{\top} + H_{\rm elastic}\idty$.

\paragraph*{Discretisation in space and time} 
The Lagrangian problem \eqref{eqn:lag-darcy-diss} is discretized in time using the first order scheme
\begin{align}
\label{eqn:num_darcy_lagrange}
\left(\kappa(\bbu^n)\frac{\bbu^{n+1}-\bbu^n}{\tau}\bJ^n,\bbw\right)=-\langle \D\bar{\calH}(\bbu^{n+1}),\bbw\rangle &&\text{for all }\bbw \in \bar{\calV}\,,
\end{align}
where $\bbu^n=\bbu(t^n,\cdot)$ is the solution at discrete time $t^n=n\tau$. This problem is discretized in space using $\bbu^n\in \bar{V}_h\subset \bar{\calV}=H^1(\bOmega,\R^d)$ via $P_k$ finite elements, where $\bar{\bs{F}}^n=\idty+\bar{\nabla}\bbu^n$ is the Lagrangian deformation gradient and $\bar{J}^n=\det\bar{\bs{F}}^n$ its determinant. The initial displacement is $\bbu^0=\bs{0}$  and we consider a constant positive $\kappa$, e.g. $\kappa(\bbu)=1$.

Similarly, we discretize the Eulerian problem \eqref{eqn:pde_eul_elastic_darcy} in time using a first order scheme 
\begin{align}
\label{eqn:num_darcy_euler}
\left(\kappa(\bu^n)\bFale^n \frac{\bu^{n+1}-\bu^n}{\tau},\bFale^n\bw\right) =& -\left\langle\D\calH(\bu^{n+1}),\bw\right\rangle   &&\text{for any }\bw \in  \calV\,.
\end{align}
where $\bs{u}^n=\bs{u}(t=n\tau,\cdot)$ is the solution at discrete time $t^n=n\tau$ and correspondingly the deformation gradient is $\bs{F}^n=(\idty-\nabla\bs{u}^n)^{-1}$. The discretisation in space is performed using $\bs{u}^n\in V_h\subset \calV=H^1(\Omega,\R^d)$ via $P_k$ finite elements. The corresponding Eulerian initial displacement is $\bu^0=\bs{0}$ and $\kappa(\bu)=1$ as above.

The difference of Lagrangian and Eulerian displacements are computed via discrete time-dependent Bochner norms $L^2(0,T;X)$
\begin{align}
\label{eqn:error_h_tau}
\mathrm{error}(h,\tau) = \left(\,\sum_{n=1}^N\tau\,\|\bu^n_h - \mathrm{proj}_{L^2}(\widehat{\bar{\bu}}^n_h)\|^2_X \right)^{1/2} ,
\end{align}
with respect to convergence in space $h\to 0$ and in time $\tau=T/N\to 0$ for a given fixed time interval $T$. Here we employ $X\in\{H^1(\Omega),L^2(\Omega)\}$ as the relevant spatial norms for the problem. With the same model parameters as in \Cref{fig:statsolconv}, the solution should converge to the same energy minimizer as $t\to\infty$.

\begin{figure}
\centering
\includegraphics[width=0.46\textwidth]{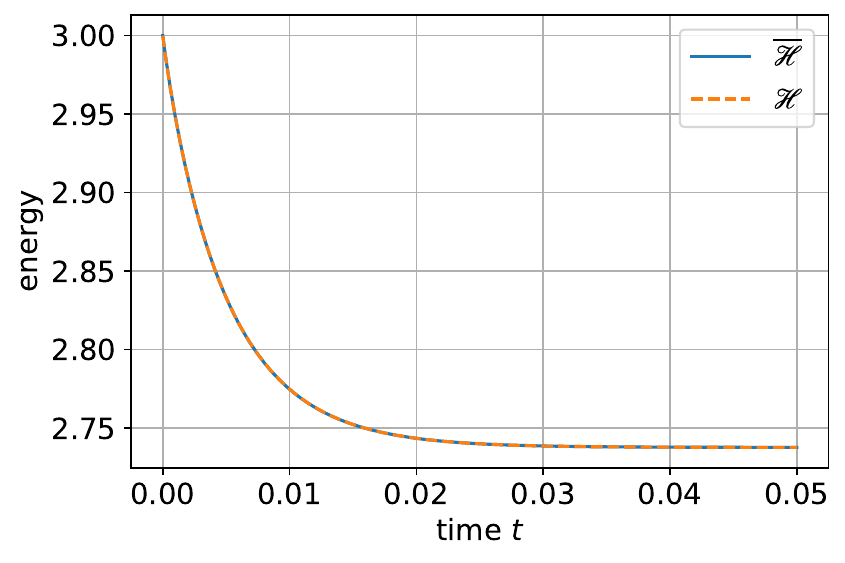}
\includegraphics[width=0.45\textwidth]{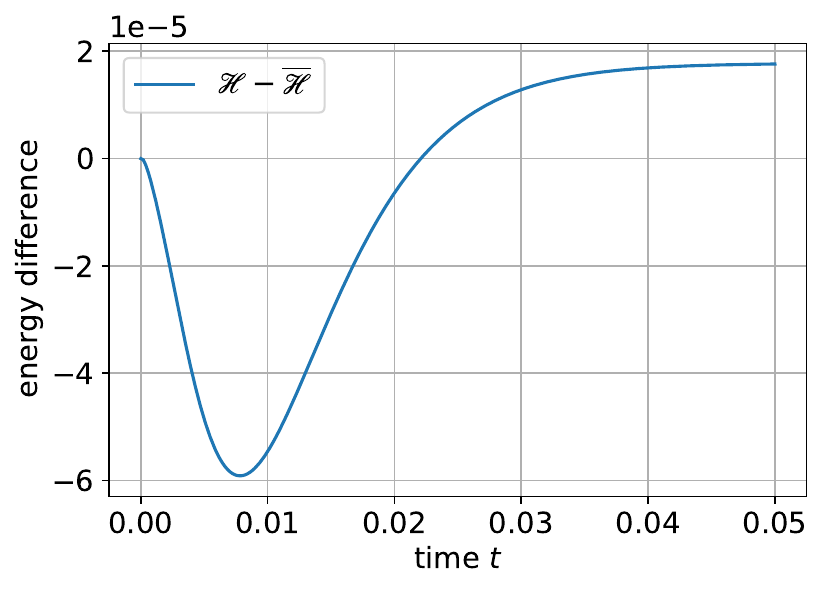}
\caption{(left) Energy $\mathscr{H}(\bq^n)$ and $\bar{\mathscr{H}}(\bbq^n)$ of numerical solutions of \eqref{eqn:num_darcy_euler} and \eqref{eqn:num_darcy_lagrange} for $P_1$ finite elements, $T=0.05$, $\tau=T/256$ on a $64\times 64$ mesh and the same parameters as in \Cref{fig:statsol} and (right) difference of energies $\mathscr{H}(\bq^n)-\bar{\mathscr{H}}(\bbq^n)$ at discrete times $t^n=n\tau$.}
\label{fig:darcy_energy}
\end{figure}

The energy descent of solutions approaching a stationary solution is shown in the left panel \Cref{fig:darcy_energy}, where at $T=0.05$ the energy is already close to its equilibrium value of $\mathscr{H}=\bar{\mathscr{H}}\approx 2.737$. The left panels shows that the difference of both solution in the energy of order $\sim 10^{-5}$, already indicating a rather good agreement.

\begin{figure}
\centering
\includegraphics[width=0.33\textwidth]{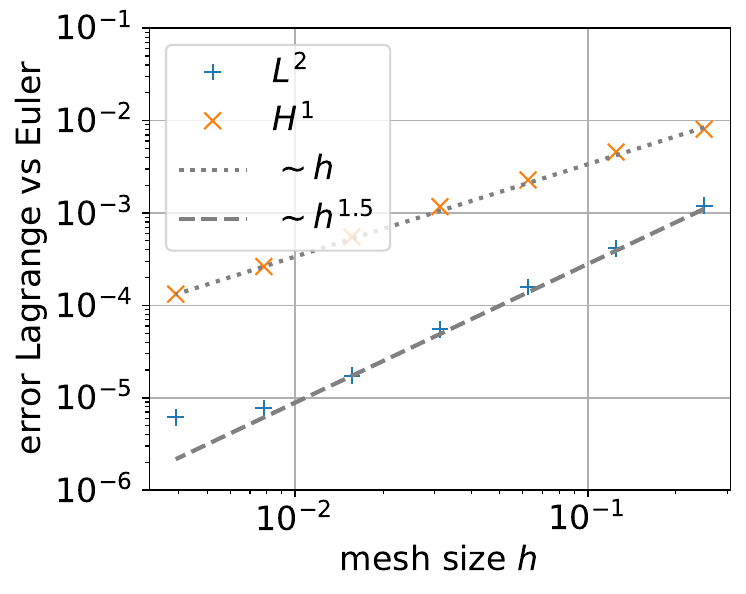}%
\includegraphics[width=0.33\textwidth]{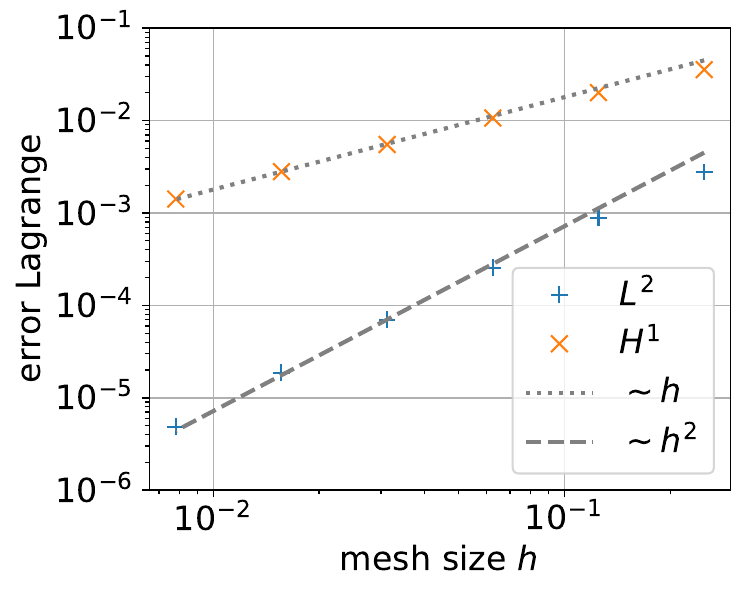}%
\includegraphics[width=0.33\textwidth]{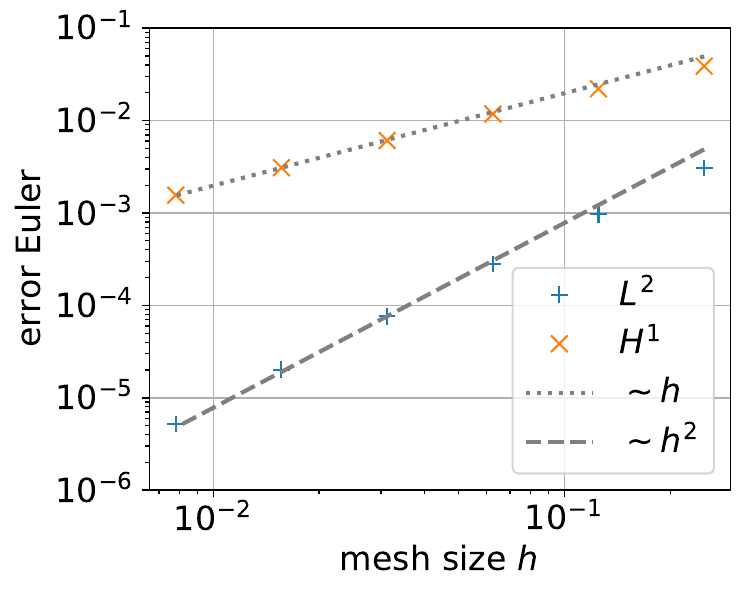}
\caption{Spatial convergence of Lagrangian and Eulerian solutions $\bbu^n$ and $\bu^n$ as $h\to 0$. (left) Difference between Lagrangian and Eulerian solutions via \eqref{eqn:error_h_tau}, (middle) corresponding Lagrangian error estimate in discrete Bochner norm $\|\bbu_{h,\tau}-\bbu_{h/2,\tau}\|_{L^2(0,T;\bar{X})}$ for $\bar{X}\in\{L^2(\bar\Omega),H^1(\bar\Omega)\}$ and (right) corresponding Eulerian error estimate in discrete Bochner norm $\|\bu_{h,\tau}-\bu_{h/2,\tau}\|_{L^2(0,T;X)}$ for ${X}\in\{L^2(\Omega),H^1(\Omega)\}$.}
\label{fig:darcy_convergence_h}
\end{figure}

We establish the experimental spatial convergence of Lagrangian and Eulerian solutions for $h\to 0$ using discrete Bochner norms
\begin{align*}
\|\bs{u}_{h_1,\tau_1}-\bs{w}_{h_2,\tau_2}\|_{L^2(0,T;X)}=\left(\sum_{n}\tau \|\bs{u}^n_{h_1,\tau_1}-\bs{w}^n_{h_2,\tau_2}\|_X^2\right)^{1/2},
\end{align*}
with time-discrete solutions $\bs{u}^n$ and $\bs{w}^n$ at a joint time $t^n=n\tau=n_1\tau_1=n_2\tau_2$, where $\bs{u}$ and $\bs{w}$ are projected onto the spatial function space of the finer uniformly refined mesh. Note that both the Lagrangian and Eulerian solutions converge optimally with order one and two in the $H^1$ and $L^2$ Bochner norms, as shown in the middle and right panel of \Cref{fig:darcy_convergence_h} respectively. The convergence of the Lagrangian to the Eulerian solution shown in the left panel is of order one in the $H^1$ norm and slightly suboptimal in the $L^2$ norm, possibly an interpolation error of the $L^2$ projection after composition with the flow map $\flowmap$. This verifies, that the implementation of the Eulerian scheme for nonlinear viscoelasticity based on the reference map provides an accurate spatial representation of the corresponding Lagrangian discretisation scheme.

\begin{figure}
\centering
\includegraphics[width=0.33\textwidth]{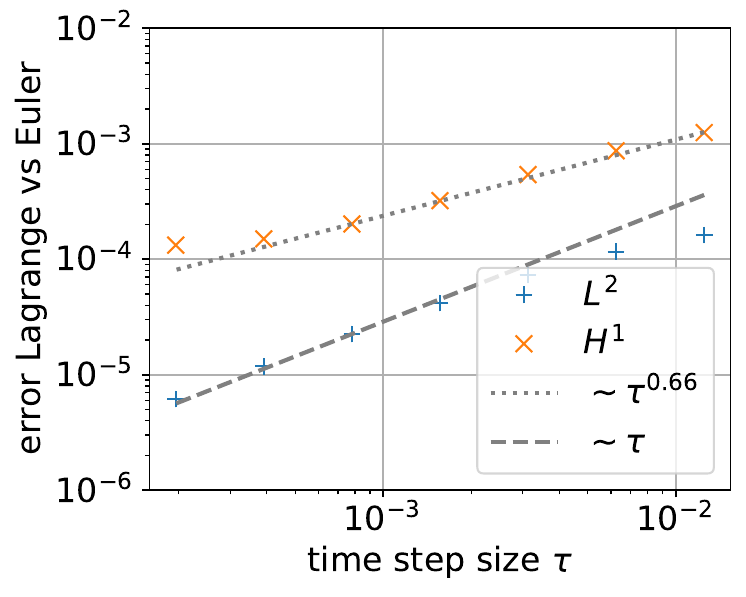}%
\includegraphics[width=0.33\textwidth]{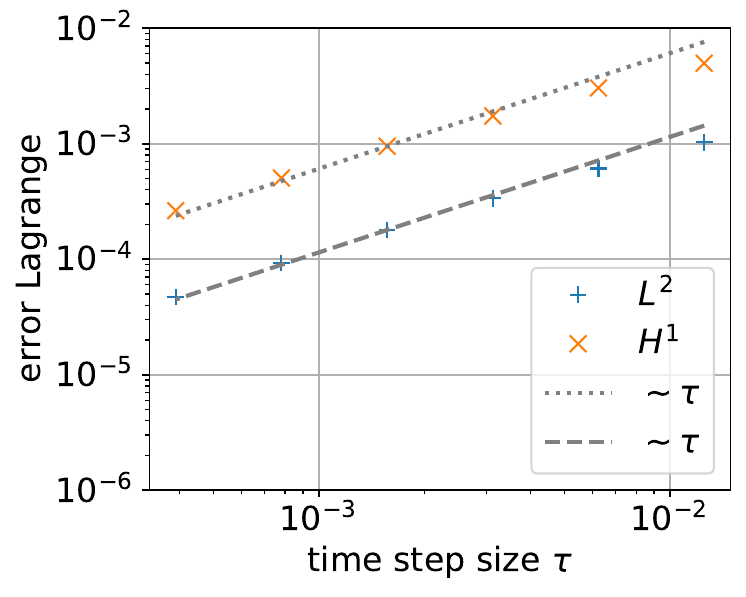}%
\includegraphics[width=0.33\textwidth]{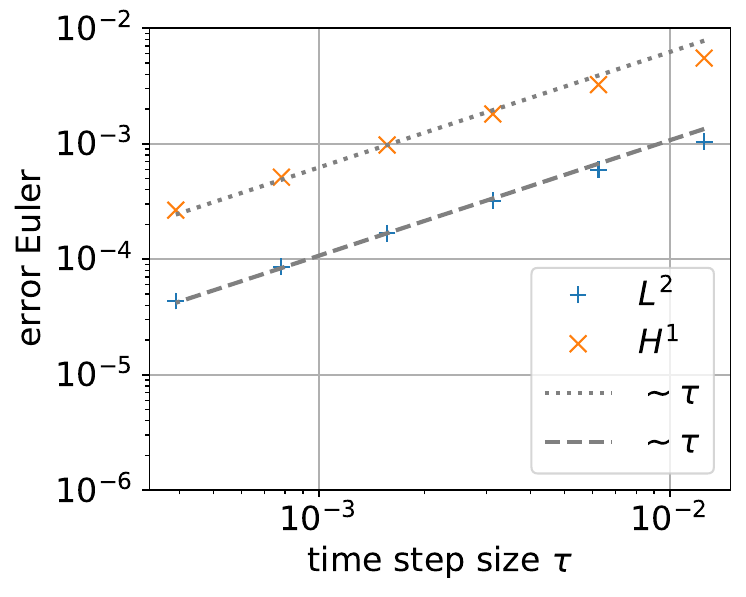}
\caption{Temporal convergence of Lagrangian and Eulerian solutions $\bbu^n$ and $\bu^n$ as $\tau\to 0$. (left) Difference between Lagrangian and Eulerian solutions via \eqref{eqn:error_h_tau}, (middle) corresponding Lagrangian error estimate in discrete Bochner norm $\|\bbu_{h,\tau}-\bbu_{h,\tau/2}\|_{L^2(0,T;\bar{X})}$ for $\bar{X}\in\{L^2(\bar\Omega),H^1(\bar\Omega)\}$ and (right) corresponding Eulerian error estimate in discrete Bochner norm $\|\bu_{h,\tau}-\bu_{h,\tau/2}\|_{L^2(0,T;X)}$ for ${X}\in\{L^2(\Omega),H^1(\Omega)\}$.}
\label{fig:darcy_convergence_tau}
\end{figure}

The temporal convergence shown in \Cref{fig:darcy_convergence_tau} for both Lagrangian and Eulerian solutions also indicate the expected first order convergence $\sim\tau$ in the Bochner norm using $L^2$ and $H^1$ spaces. This order now slightly deteriorates for the convergence via \eqref{eqn:error_h_tau} in the $H^1$ norm, possibly again due to an interpolation error due to the $L^2$ projection $\mathrm{proj}_{L^2}$ of Lagrangian functions to their corresponding Eulerian representation. This shows that the gradient flow based discretisation based on incremental time minimisation also provides an accurate temporal representation of the Lagrangian discretisation scheme and the corresponding solution.

Let us make a short remark on how such a gradient flow system for a viscoelastic system might emerge from a full hyperelastic material with inertia via an overdamped limit.

\begin{remark}[Gradient systems from overdamped limit of damped Hamiltonian systems]\label{rem:overdamped-ham}
    Gradient systems as \eqref{eqn:GS_equation}, that have a flow map, a deformation or a displacement within their variables can sometimes be regarded as limit of overdamped Hamiltonian systems. To fix the ideas, consider the structure in \eqref{eqn:MAU} and let us introduce a new skew-symmetric bilinear form $\bar{j}(\bar{\bs{\eta}}, \bar{\bw}):\bcalW\times\bcalW\to\R$, which is usually associated to a Poisson operator, see e.g. \cite{peschka2022variational}. The evolution equation for a damped Hamiltonian system can then be written in the saddle point structure 
    \begin{subequations} \label{eqn:ext-damp-ham1}
        \begin{align}
            \bar{k}(\bar{\bs{\eta}},\bar{\bs{\xi}}) - \bar{j}(\bar{\bs{\eta}},\bar{\bs{\xi}}) + \bar{b}(\dot{\bbq},\bar{\bs{\xi}}) &= 0\,, &&\text{ for all } \bar{\bs{\xi}}\in\bcalW
            \\
            \bar{b}(\bbw,\bar{\bs{\eta}})
            &=\left\langle\D\bcalH(\bar{\bq}),\bbw\right\rangle\,  &&\text{ for all } \bar{\bw} \in \bcalV\,.
        \end{align}
    \end{subequations}
    As an example, consider the vector of state variables $\bbq=(\flowmap, \bmom)$, where $\flowmap$ is a flow map as defined in subsection \ref{sec:kinematics} and $\bmom:[0,T]\times\bOmega\to\R^d$ is the momentum density. We split the contributions of the Hamiltonian energy density into kinetic and internal energy:
    \begin{subequations}
    \begin{equation*}
        \bar{\calH}(\bbq):=\int_{\bOmega}\frac{\vert\bmom\vert^2}{2\brho} + \bar{W}(\bbF)\,\dbbx\,,
    \end{equation*}
    with $\brho$ being the mass density and $\bbF=\bnabla\flowmap$ the deformation gradient as before.
    While the bilinear forms are so defined:
    \begin{align*}
        \bar{b}(\bbv,\bbw):=& \int_{\bOmega} \bbv_{\flowmap}\cdot\bbw_{\flowmap} + \bbv_{\bmom}\cdot\bbw_{\bmom}\,\dbbx\,,
        \\
        \bar{j}(\bbw_1, \bbw_2):=& \int_{\bOmega}  \bbw_{1,\bmom}\cdot\bbw_{2,\flowmap} -\bbw_{1,\flowmap}\cdot\bbw_{2,\bmom} \,\dbbx\,,
        \\
        \bar{k}(\bbw_1, \bbw_2):=& \int_{\bOmega} \mu(\bbx)\bnabla\bbw_{1,\bmom}\cdot\bnabla\bbw_{2,\bmom}\,\dbbx\,,
    \end{align*}
    \end{subequations}
    Altogether these ingredients and \eqref{eqn:ext-damp-ham1} give a weak formulation of a damped Hamiltonian system for a viscoelastic material in Lagrangian coordinates.
    If we are interested in the long-time behaviour of highly viscous materials, then it is possible to rescale time by $t = \hat{t}/\varepsilon$ with $\varepsilon \ll 1$ and subsequently the viscosity coefficient $\mu$, see e.g. \cite{zafferi2022generic}. We then obtain a PDE system where $\dot{\bmom}$ can be neglected. The resulting equations can be expressed as a \emph{overdamped} (od) gradient system solely in terms of the variable $\bbq_{\rm od}=\flowmap$. In that case, the free energy simplifies to 
    \begin{subequations}\label{eqn:od_setting}
        \begin{align}
        \bcalH_{\rm od}(\bbq_{\rm od}):= \int_{\bOmega}W(\bbF)\,\dbbx\,,
    \end{align}
    while with the space $\bcalW=0$, the dissipation is completely encoded in the bilinear form 
    \begin{equation}
        \bar{s}_{\rm od}(\bbv_1,\bbv_2):= \bar{k}(\bbv_1,\bbv_2)\,.
    \end{equation}
    \end{subequations}
    Equations \eqref{eqn:od_setting} suffice to obtain a weak formulation of a gradient system as in \eqref{eqn:MAUS}.
    This example was introduced to highlight that there might be cases where similar processes can be encoded into different bilinear forms, perhaps after a formal (overdamped) limiting procedure, and to present a case where the transformation mechanism of Section \ref{sec:trafo} does not apply.
\end{remark}


\section{Geophysical flows with reference map and diffusion}
\label{sec:extension_concentration}
In this section we expand the previous model to classes of models that incorporate various geophysical effects and contain a concentration variable, that models the hydration of rock, \emph{i.e.}, the transport of a fluid content. We will thus extend the energy and the dissipative structure previously presented and focus more on the Eulerian gradient structure that emerges from a transformation of a Lagrangian structure. 
The goal of  \Cref{sec:extension_concentration} is to move towards more complex and geophysically justified models, while preserving the thermodynamic gradient flow structure.

\subsection{Eulerian poroelasticity}
\label{sec:eul-poro}
Firstly, we want to construct a gradient flow formulation of a nonlinear poroelastic model in the saddle point format in the Eulerian frame by mapping from a Lagrangian frame. Therefore, we are not going to reiterate the complete Lagrangian and Eulerian model but focus on the crucial steps that need to be taken in addition to the computations that we already made in the preceding section. 

\subsubsection*{Eulerian formulation by transformation}

Firstly, assume we have a Lagrangian formulation in terms of displacement and (molar) concentration, \emph{i.e.}, $\bbq=(\bbu,\bc)\in\bcalQ$ with displacement $\bbu:\bOmega\to\R^d$ and concentration $\bc:\bOmega\to\R^n$ for $\bc=(c_1,...,c_n)$ for each component $c_i$ for $i=1,...,n$.
The system is described by a corresponding Hamiltonian functional  $\bar{\calH}:\bcalQ\to\R$. In correspondence to \eqref{eqn:MAUS}, we assume $\bbw=\bbw_{\bc}:\bOmega\to\R^n\in\bcalW$,  so that we seek $\bbq:[0,T]\to\bcalQ$ and $\bbeta:[0,T]\to\bcalW$ such that
\begin{align*}
        \bar{k}(\bbeta,\bbw) + \bar{b}(\dot{\bbq},\bbw) &= 0 &&\forall\, \bbw\in\bcalW\,,
        \\ 
        \bar{b}(\bbv,\bbeta) - \bar{s}(\dot{\bbq},\bbv) &=\left\langle\D\bcalH(\bbq),\bbv\right\rangle_{\bcalV}  &&\forall\,\bbv \in \bcalV\,.
\end{align*}
where as before $\dot{\bbq}=(\dot{\bbu},\dot{\bc})=(\partial_t\bbu,\partial_t\bc)$. We fix the form of $\bar{b}$ to be 
\begin{align*}
\bar{b}(\dot{\bbq},\bbw):=& \int_{\bOmega} \dot{\bc}\cdot\bbw\,\dbbx\,,
\end{align*}
We set $\bq=(\bu, \con)\in\calQ$ as state vector, defined by the transformation $\trafo:\bcalQ\to\calQ$ via
\begin{equation}\label{eqn:def-trafo-el-di}
    \buale(t,\flowmap(t,\bbx)):= \bbu(t,\bx)\in\R^d\,, \qquad\con(t, \flowmap(t,\bbx)) := \frac{\bc(t,\bbx)}{\bJ(t,\bbx)}\in\R^n\,.
\end{equation}
This form of the mapping, specifically for the concentration, is motivated by the molar concentration being an extensive variable that then needs to be mapped like this and leads to a continuity equation, as we discussed extensively in \cite{zafferi2022generic}.
The free energy of the system will have to satisfy the closure condition $\calH(\trafo(\bbq))=\bar{\calH}(\bbq)$.  We write the weak formulation in terms of functions $\bv=(\bv_{\bu}, \bv_{\con})\in\calV$ and with the auxiliary space $\bw=\bw_{\con}\in{\calW}$.
As we have outlined in \Cref{sec:trafo} the main element needed to identify the bilinear forms $\brm, \krm $ and $\srm$ is the linearisation $\dtrafo:\bcalV\to\calV$ of the transformation $\trafo:\bar{\calQ}\to\calQ$ from \eqref{eqn:def-trafo-el-di},
\begin{align}
    \dtrafo(\bbq):=\begin{pmatrix}
        \nabla\ale\Box & 0 \\
        -\nabla\cdot(\con\Box) & \frac{1}{\jale}
    \end{pmatrix}\,.
\end{align}
The similarly mapping $\bs{\eta}(t,\flowmap(t,\bbx)):= \bbeta(t,\bx)$ defines a linearisation $\dtrafo_{\calW}:\bcalW\to\calW$ via $\dtrafo_{\calW}(\bbq)\bbeta=\bs{\eta}$ or simply $\dtrafo_{\calW}=1$. This allows us to define the mapped bilinear form by the identity
\begin{align}
\label{eqn:trafo_identity}
\brm(\dtrafo \bbv,\dtrafo_{\calW} \bbw) := \bar{b}(\bbv,\bbw),
\end{align}
which produces with $\bw\in\calW$ and $\bv=(\bv_\bu,\bv_c)\in\calV$
\begin{align}\label{eqn:eul-b}
\brm(\bv,\bw)=\int_\Omega \bw\cdot(\bv_c + \nabla\cdot(\con\bF\bv_{\bu}))\,\dbx\,,
\end{align}
because then we satisfy the transformation identity \eqref{eqn:trafo_identity} via 
\begin{align*}
\begin{split}
\brm(\dtrafo \bbv,\dtrafo_{\calW} \bbw) &= \brm\Bigl(\bigl((\nabla\ale)\bbv_u,-\nabla\cdot(\con\bbv_u)+J^{-1}\bbv_c\bigr),\bbw\Bigr) \\
&=\int_\Omega \bbw\cdot(J^{-1}\bbv_c \underbrace{- \nabla\cdot(\con\bbv_u) + \nabla\cdot(\con\overbrace{\bF(\nabla\ale)}^{=\idty}\bbv_u}_{=0})\,\dbx\\
&=\int_{\bOmega} \bbw\cdot\bbv_c\,\dbbx = \bar{b}(\bbv,\bbw).
\end{split}
\end{align*}
This is a key observation, as the bilinear form $\brm$ generates convective derivatives and coupling of the equation of displacement and concentration. Similarly, the other bilinear forms are defined by
\begin{align*}
\krm(\dtrafo_{\calW} \bbw_1,\dtrafo_{\calW} \bbw_2) := \bar{k}(\bbw_1,\bbw_2), \qquad 
\srm(\dtrafo \bbv_1,\dtrafo \bbv_2) := \bar{s}(\bbv_1,\bbv_2).
\end{align*}
\subsubsection*{Weak and strong formulation}
With these identities relating the Lagrangian and the Eulerian formulation, we can state the entire weak formulation of the system coupling elasticity and diffusion, where we seek $\bq:[0,T]\to\calQ$ with $\partial_t\bq\in\calV$ solving
\begin{subequations}
\label{eqn:FULL_MAUS}
\begin{align}
        \krm(\bs{\eta},\bw) + \brm(\partial_t\bq,\bw) &= 0 &&\forall\, \bw\in\calW\,,
        \\ 
        \brm(\bv,\bs{\eta}) - \srm(\partial_t\bq,\bv) &=\left\langle\D\calH(\bq),\bv\right\rangle_{\calV}  &&\forall\,\bv \in \calV\,.
\end{align}
\end{subequations}
The bilinear forms that we use above in \eqref{eqn:FULL_MAUS} are defined
\begin{subequations}
\label{eqn:FULL_BILINEAR}
\begin{align}
    \label{eqn:diff_kblf}&\krm(\bw_1,\bw_2):=
    \int_{\Omega}\,\nabla\bw_{1}\cdot D(\con)\nabla\bw_{2}\,\dbx\,,
    \\
    \label{eqn:darcy-sblf}&\srm(\bv_1,\bv_2):= \int_{\Omega}\nu\nabla_s(\bF\bv_{1,\bu})\cdot\nabla_s(\bF\bv_{2,\bu})\,\dbx\,,
    \\
    \label{eqn:diff_bblf}&\brm(\bv,\bw):=\int_\Omega \bw\cdot(\bv_c + \nabla\cdot(\con\bF\bv_{\bu}))\,\dbx\,,
\end{align}
\end{subequations}
where we used $\bv_i=(\bv_{i,\bu},\bv_{i,c})\in\calV$ and $\bw_i\in\calW$ for $i\in\{1,2\}$ and the symmetric gradient $\nabla_s \bs{a}=\tfrac12 \nabla\bs{a}+\tfrac12(\nabla\bs{a})^T$. The diffusion tensor is $D(\con)\in\R^{n\times n}$ is symmetric positive semidefinite. Note that the bilinear forms are linear in $\bv_i,\bw_i,\bv,\bw$ but also have a general dependence on $\bq$ through $\con$ and $\bF=(\idty-\nabla\bu)^{-1}$. 
For the free energy of the system we choose
\begin{subequations}\label{eqn:hamiltonian-eul-conc}
\begin{align}\label{eqn:energy_single_phase}
\begin{split}
    &\calH(\bq):=\int_{\Omega} H_{\rm elastic}(\bF) + {H}_{\rm Biot}(\con,J) + {H}_{\rm gravity}(\con) + \frac{\varepsilon}{2}|\nabla\con|^2\,{\rm d}{\bx}\,,
    \\
    &\text{where}\quad H_{\rm elastic}(\bF):= \frac{\mu}{2}\,{\rm tr}(\bC_\textrm{iso} - \idty)\,,
    \quad 
    H_{\rm Biot}(\con,J):= \frac12\,\pi^2\,,
    \quad 
    H_{\rm gravity}(\con):=g \bs{m}\cdot\con z\,,
\end{split}
\end{align}
with coordinate vector $\bx=(x,z)\in\R^d$ where $x\in\R^{d-1}$ and $z\in\R$, $\bs{g} = (0,g)$ where $0\in\R^{d-1}$ and $g\in\R$,  deformation gradient $\bF=(\nabla\ale)^{-1}=(\idty-\nabla\bu)^{-1}$, the determinant of the deformation gradient $J=\det\bF$, and the isochoric right Cauchy-Green deformation tensor $\bC_\textrm{iso}=\bF^T\bF/J^{2/d}$,  with pore pressure $\pi=\bs{\alpha}\cdot(\con-\con^*)-(J-1)$, and $\bs{m}=(m_1,...,m_n)\in\R^n$ the vector of molar masses and $\bs{\alpha}=(\alpha_1,...,\alpha_n)\in\R^n$ the vector of molar volumes.
The energy functional is composed of four distinct terms: $H_{\rm elastic}$ incorporates the elastic material response through the isochoric part of the right Cauchy-Green deformation tensor $\bC_\textrm{iso}$, $H_{\rm Biot}$ encodes the coupling of displacement and concentration through volume changes $\jale$ with concentration variations from the steady state $\con^*$ by a Biot pressure $\pi$, $H_{\rm gravity}$ encodes gravitational forces in the direction of the $d$-th coordinate and we include a higher-order regularising term for $\con$.
\end{subequations} 
The total energy \eqref{eqn:hamiltonian-eul-conc} generates the following driving force 
\begin{align*}
\begin{split}
    \left\langle\D\calH(\bq), \bv\right\rangle &= \int_{\Omega}\bF^{\top}\partial_{\bF}\left(H_{\rm elastic}\right)\bF^{\top}{:}\nabla\bv_{\bu} - \pi\jale\bF^{\top}:\nabla\bv_{\bu} + \pi \bs{\alpha}\cdot \bv_{\con}+ g z \bs{m}\cdot \bv_{\con} + \varepsilon \nabla\con\cdot\nabla \bv_{\con}\dbx
    \\
    &= \int_{\Omega}\Big(\tfrac{\mu}{2}\left(2\bC_\textrm{iso}-\tfrac{2}{d}{\rm tr}(\bC_\textrm{iso})\idty \right)- \pi\jale\Big)\bF^{\top}{:}\nabla\bv_{\bu} 
    + (\pi\bs{\alpha} + g z \bs{m} )\cdot \bv_{\con} + \varepsilon \nabla\con\cdot\nabla \bv_{\con}\dbx
\end{split}
\end{align*}
With the definitions above we can state the entire weak formulation \eqref{eqn:FULL_MAUS} for the separate components of the solution and test functions as follows:
\begin{subequations}
\label{eqn:full_weak}
\begin{align}
\int_\Omega \bs{\eta}\nabla\cdot(\con\bF\bv_{\bu}) -\nu\nabla_s(\bF\partial_t \bu)\cdot(\nabla_s\bF\bv_{\bu})\,\dbx  = \int_\Omega \left(\bF^{\top}\partial_{\bF}H_{\rm elastic} - \pi\jale\idty\right)\bF^{\top}:\nabla\bv_{\bu} \dbx,\\
\int_\Omega \bw\Big(\partial_t\con + \nabla\cdot(\con\bF\partial_t\bu)\Big)\dbx = -\int_\Omega D(c)\nabla\bs{\eta}\cdot\nabla\bw\,\dbx,\\
\int_\Omega \bs{\eta}\bv_{\con}\dbx = \int_\Omega (\bs{\alpha}\pi + \bs{m}gz)\cdot\bv_{\con} + \varepsilon\nabla\con\cdot\nabla\bv_{\con}\,\dbx.
\end{align}
\end{subequations}
which is complemented with suitable boundary conditions on $\partial\Omega$. The strong form of the system \eqref{eqn:full_weak} yields the following model for a poroelastic material
\begin{subequations}
\label{eqn:strong}
\begin{align}
    \label{eqn:strong-u}
    - \nu\,{\rm div}\nabla_{\rm s}(\bF\dot{\bu}) + \con\cdot\nabla{\bs{\eta}} =& {\rm div}(\bs{\sigma})\,,
        \\\label{eqn:strong-c}
        \dot{\con}+ \mathrm{div}(\con \bF \dot{\bu}) =& {\rm div}\left(D(\con)\nabla{\bs{\eta}} \right)\,,
\end{align}
\end{subequations}
where we have the elastic Cauchy stress
$\bs{\sigma}=\partial_{\bF}H_{\rm elastic}\bF^{\top} + H_{\rm elastic}\idty + ( \tfrac{1}{2}\pi^2-\jale\pi)\idty$, 
the pore pressure 
${\pi}=\bs{\alpha}\cdot(\con - \con^*) - (\jale - 1)$ and the chemical potential 
$\bs{\eta} = \pi\bs{\alpha} + g \bs{m} z - \varepsilon\Delta\con$.
This is an Eulerian model supporting large deformations and it models the viscous relaxation of the compressible elastic body by shear forces using a Biot-type pore pressure encoded in the Cauchy stress $\bs{\sigma}$. Remember that in this model the Eulerian solid velocity is $\bv=\bF\dot{\bu}$.
\begin{remark}
By defining the Lagrangian elastic energy by $\bar{H}_{\rm elastic} = H_{\rm elastic}\jale$ one recovers the classical elastic Cauchy stress $\bs{\sigma}$ in \eqref{eqn:strong-u} and, in a similar manner, the contribution to the thermodynamic pressure due to ${H}_{\rm Biot}$. 
One essential step in the derivation of the Cauchy stress is the following computation, where for simplicity we abbreviate $H = H_{\rm elastic}$:
\begin{align*}
    \bF^{-\top}{\rm div}\left(\bF^{\top}\frac{\partial H}{\partial\bF}\bF^{\top}\right) &= \bF^{-\top}_{ij}\nabla_k\left(\bF^{\top}_{jl}\frac{\partial H}{\partial\bF_{lm}}\bF^{\top}_{mk}\right) 
    \\ &= \bF^{-\top}_{ij}\bF^{\top}_{jl}\nabla_k\left(\frac{\partial H}{\partial\bF_{lm}}\bF^{\top}_{mk}\right) + \bF^{-1}_{ji}\frac{\partial H}{\partial\bF_{lm}}\bF_{km}\bar{\nabla}_j(\bF_{ln})\bF^{-1}_{nk} 
    \\ &= {\rm div}\left(\partial_{\bF}H\bF^{\top}\right) + \frac{\partial H}{\partial\bF_{lm}}\bar{\nabla}_j(\bF_{lm})\bF^{-1}_{ji}
    = {\rm div}\left(\partial_{\bF}H\bF^{\top} + H\idty\right)\,,
\end{align*}
where the sum over the same indices has to be taken, an identity relating Eulerian spatial derivatives with the third order deformation tensor $\nabla_k\bF_{lj} = \bar{\nabla}_{s}\bF_{lj}\bF^{-1}_{sk},\, \bar{\nabla}_s\bF_{lj}=\bar{\nabla}^2_{sj}\flowmap_l$ with $\flowmap$ being the flow map was used, and, in the last equality, we applied the chain rule for $\nabla H$. 
\end{remark}

This is a higher-order diffusive model, in which the role of the fourth-order term proportional to $\varepsilon$ is to regularize and stabilize effects due to gravity. The mobility $D(\bs{c})$ is usually degenerate and diagonal, \emph{i.e.}, $D_{ij}(\bs{c})=|c_i|^m D^0_i \delta_{ij}$ for some exponent $m\ge 0$, $D_i^0>0$ and $i,j=1,...,n$. For ideally (diffusing) gases, often logarithmic terms in the energy can be used to enforce positivity of the components $c_i$. In the following examples we restrict to a single component $n=1$ and two spatial dimensions $d=2$.

\subsection{Discretisation}
We discretise the poroelastic model using a variational, structure-preserving approach, \emph{i.e.}, we start the discretisation not at the equation level but starting from the general saddle-point structure \eqref{eqn:FULL_MAUS}. Here we emphasize that in addition, the bilinear forms $\krm,\srm,\brm$ also depend explicitly on the state $\bq=(\bu,\con)$, a dependence that we omitted for brevity. However, for the time-discretisation we are going to explicitly highlight the state-dependence by writing $\krm(\bq;\bw_1,\bw_2)$, $\srm(\bq;\bv_1,\bv_2)$ and $\brm(\bq;\bv,\bw)$. Furthermore, we assume we are given a sequence of time-steps $0=t^0<t^1<...<t^N=T$ and write $\bq^k(\bx)=\bq(t^k,\bx)$ and $\bs{\eta}^k(\bx)=\bs{\eta}(t^k,\bx)$ or also simply $\bq^k=\bq(t^k)$ and $\bs{\eta}^k=\bs{\eta}(t^k)$. We approximate the time-derivative $\partial_t\bq(t^k)=(\bq^{k}-\bq^{k-1})/\tau^k$ with $\tau^k=t^{k}-t^{k-1}$.

We discretise \eqref{eqn:FULL_MAUS} using finite elements in space and semi-implicitly in time, \emph{i.e.}, the derivative of the energy $\D\calH(\bq)$ is discretised implicitly and the state-dependence of the bilinear forms in \eqref{eqn:FULL_BILINEAR} is discretised explicitly in $\bq$ and implicitly in the remaining variables. Therefore, for a given triangulation $\mathcal{T}_h$ of the domain $\Omega=\bigcup_{T\in\mathcal{T}_h}T$ introduce the generic finite-dimensional function space 
\begin{align}\label{eqn:FEM-space}
V_{h}^{m,k}=\big\{\bs{\varphi}=(\varphi_1,...,\varphi_m)\in C^0(\Omega;\mathbb{R}^m)\,:\,\varphi_i|_T \in P_k(T),\, T\in \mathcal{T}_h,\,i\in\{1,...,m\}\big\}\,,
\end{align}
where $m=1$ for scalar functions and $m=d$ for vector fields and $k$ is the polynomial degree of the $H^1$ conforming function on elements $T\in\mathcal{T}_h$. We use the spaces 
\begin{align}\label{eqn:approx-state}
\bq^k = (\bs{u}^k,\bs{c}^k) = V_{h}^{d,2} \times V_{h}^{1,1} =: \mathcal{V}_h \qquad\text{and}\qquad \bs{\eta}^k = V_{h}^{1,1} =: \mathcal{W}_h\,,
\end{align}
where the space for $\bs{u}^k$ is complemented with suitable Dirichlet boundary conditions $\bs{u}^k\cdot\bs{\nu}=0$ (sliding) or $\bs{u}^k=0$ (no-slip) to be specified later on parts of $\partial\Omega$. For a given $\bs{q}^{k-1}$, we seek $\bq^k$ and $\bs{\eta}^k$ such that
\begin{subequations}  \label{eqn:DMAUS}
    \begin{align}\label{eqn:DMAUS-1}
        \krm(\bq^{k-1};\bs{\eta}^k,\bw) &+ \brm(\bq^{k-1};\tfrac{\bq^k-\bq^{k-1}}{\tau^k},\bw) &&= 0 &&&&\text{for all } \bw\in\calW_h\,,
        \\ \label{eqn:DMAUS-2}
        \brm(\bq^{k-1};\bv,\bs{\eta}^k) &- \srm(\bq^{k-1};\tfrac{\bq^k-\bq^{k-1}}{\tau^k},\bv) &&=\left\langle\D\calH(\bq^k),\bv\right\rangle_{\calV}  &&&&\text{for all } \bv \in \calV_h\,.
    \end{align}
\end{subequations}
In a previous work \cite{schmeller2023gradient} we showed that such discretisation scheme implements an incremental minimisation scheme for $\bq^k$, where the energy is provable decreasing. The problem \eqref{eqn:DMAUS} results in a nonlinear saddle-point problem that is solved using a Newton method. Here, we solve the problem using a constant constant time step size $\tau^k\equiv \tau$ for all $k\in\{1,...,N\}$. \Cref{eqn:DMAUS} is a semi-implicit \emph{monolithic structure-preserving discretisation} of the gradient flow equation \eqref{eqn:Onsager_equation} or the corresponding saddle-point structure \eqref{eqn:FULL_MAUS}. Next we are going to present some exemplary solutions of this numerical scheme that are relevant in the context of geophysical applications.

\subsection{Porosity waves}\label{sec:pw-single-phase}
\begin{figure}[b!]
\centering
\begin{minipage}{0.49\textwidth}
\includegraphics[width=\textwidth]
{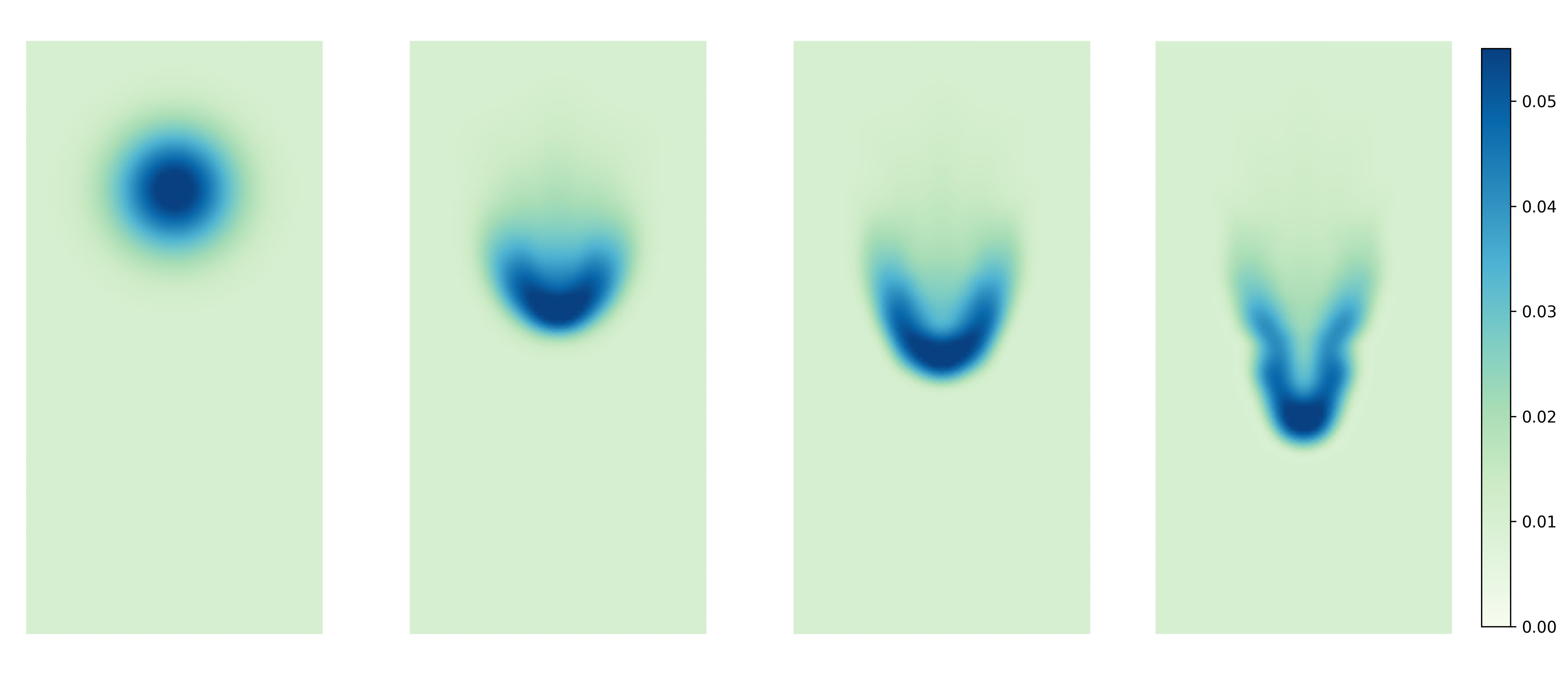}
\end{minipage}\hfill
\begin{minipage}{0.49\textwidth}
\includegraphics[width=\textwidth]
{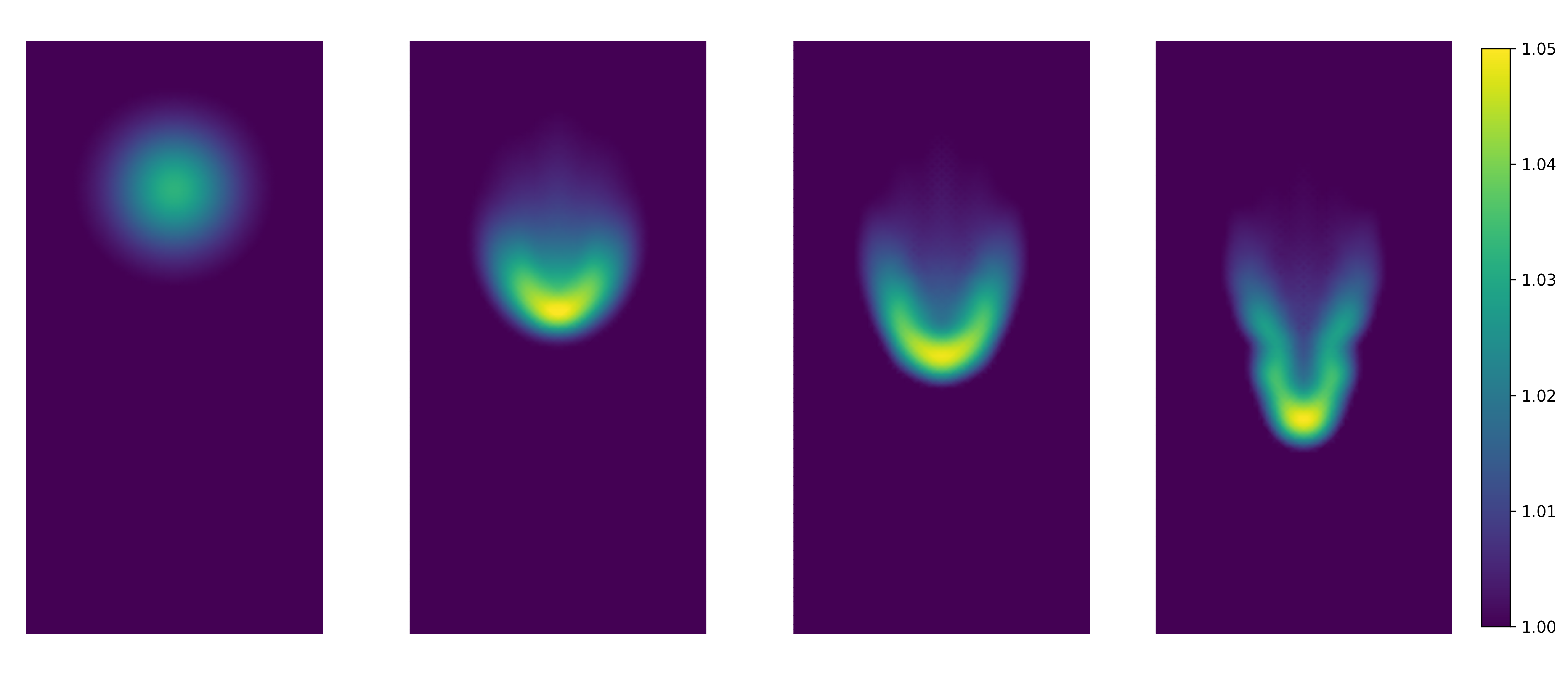}
\end{minipage}
\caption{Porosity wave moving fluid downwards and displacing the (elastic)  material in the \textbf{flowing regime} at times $t=0,10/3,20/3,40/3$. On the left, Concentration $\con$ and on the right the determinant $J=\det\bF$. Both for time advancing from left to right.}
\label{fig:flowing}
\end{figure}

While the model \eqref{eqn:strong}, its weak form and the corresponding discretisation appear straightforward, there are a couple of peculiarities that we briefly want to comment on. Firstly, note that the diffusion equation is of fourth order due to the regularisation term in the energy, which also contributes a third-order Korteweg-like term to the momentum balance via the chemical potential $\bs{\eta}$. Both the fourth-order term in the diffusion and the third-order term in the momentum equation are consistently handled via a mixed finite-element scheme by introducing the chemical potential as a separate variable. Using this variational structure, we are able to introduce additional physics, e.g., gravitational forces, in the momentum and diffusion equation leading to the correct equilibrium states. For diffusion with gravity, the corresponding equilibrium is given by the barometric formula. However, note that while $\bs{\eta}$ is the correct chemical potential for this thermodynamic system, $\pi$ is not the full thermodynamic pressure but only a Biot-type contribution to it.

In the following we are going to show some typical solutions for excess concentration traversing through a poroelastic material driven by gravity. 
We solve this problem in two dimensions in a rectangular domain $\Omega=(0,L)\times(0,H)$ with components $\bs{x}=(x_1,x_2)\in\Omega$ and with initial data
\begin{align}\label{eqn:inital_datum1}
\con_0({x})= \bar{c}_0 + \bar{c}_1 \exp(-\tfrac{1}{r^2}|\bs{x}-\bs{x}_0|^2)\,,
\end{align}
and an initial displacement $\bs{u}_0(\bs{x})$ that minimises the energy 
\begin{align*}
\calH_0(\bs{u}_0):=\calH\bigl((\bs{u}_0,\con_0)\bigr) + \int_\Omega |\nabla\bs{u}_0|^2\,\mathrm{d}x\,,
\end{align*}
for the given initial concentration $\con_0$. 
We employ slip boundary conditions $\bs{u}(t)\cdot\bs{\nu}=0$ on the side walls $x_1\in\{0,L\}$ and no-slip boundary conditions on the top and bottom walls $x_2\in\{0,H\}$. In order to conserve fluid content $\con$, we employ homogeneous natural boundary conditions. The nondimensional domain size is $L=1$ and $H=2$ and in the initial data we use $r=1/5$. We use $\bar{c}_0=10^{-2}$ and $\bar{c}_1=5\cdot 10^{-2}$ and viscosity $\nu=10^{-3}$.
For gravity we use $\bs{g}=(0,g)$, where $g$ and the remaining parameters are given in \Cref{tab:regime_comparison}. Both simulations use 300 time steps to reach the final time $T$.
\begin{table}[htbp]
\centering
\caption{Comparison of selected parameters for flowing and diffusive regimes.}
\begin{tabular}{c|rr}
\hline
\textbf{Parameter} & {flowing regime} & {diffusive regime} \\
\hline
$g$                      & $0.1$           & $0.01$ \\
$D_0$   & $0.1$            & $1000$ \\
$\varepsilon$             & $10^{-3}$        & $10^{-5}$ \\
$\mu$                     & $10^{-9}$        & $10^{-3}$ \\
$T$ & $10.0$ & $4.0$ \\
\hline
\end{tabular}
\label{tab:regime_comparison}
\end{table}

\begin{figure}[hb!]
\centering
\begin{minipage}{0.49\textwidth}
\includegraphics[width=\textwidth]
{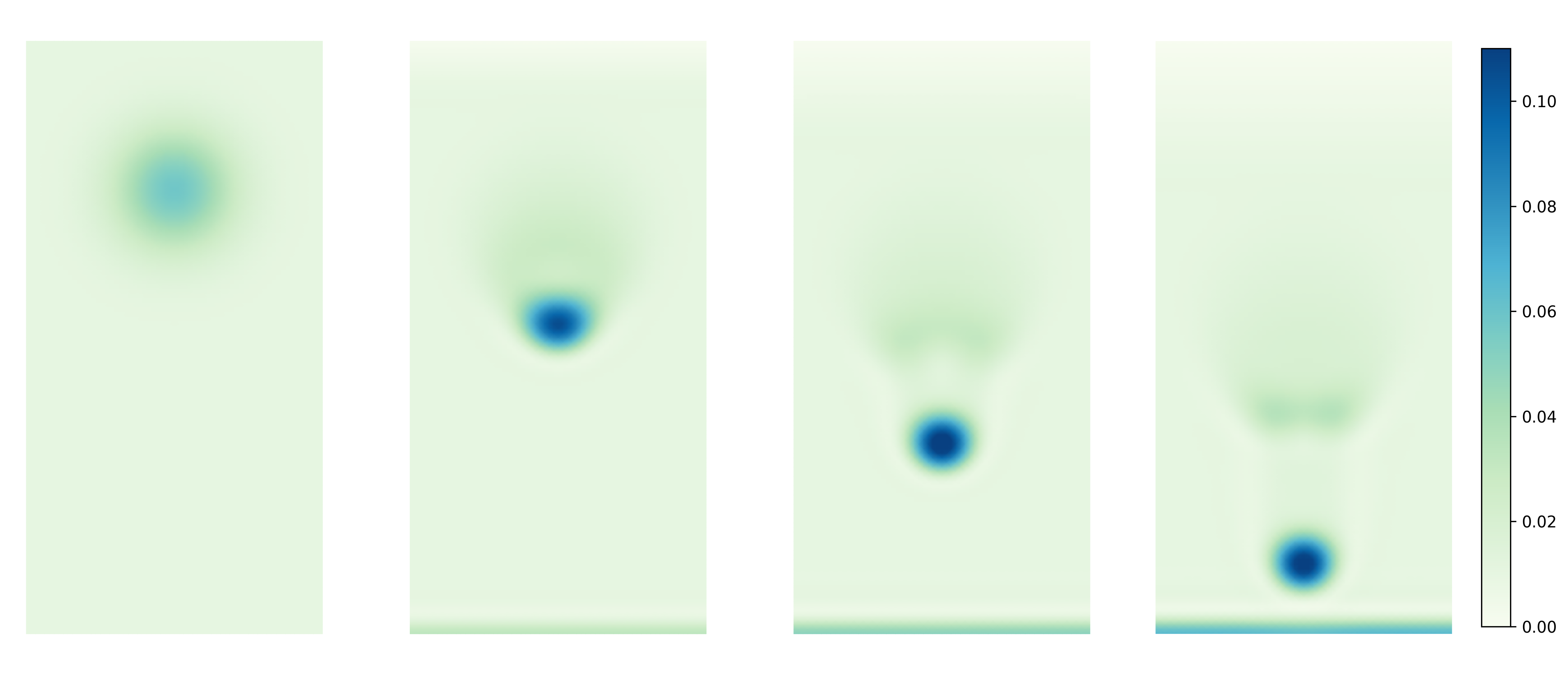}
\end{minipage}\hfill
\begin{minipage}{0.49\textwidth}
\includegraphics[width=\textwidth]
{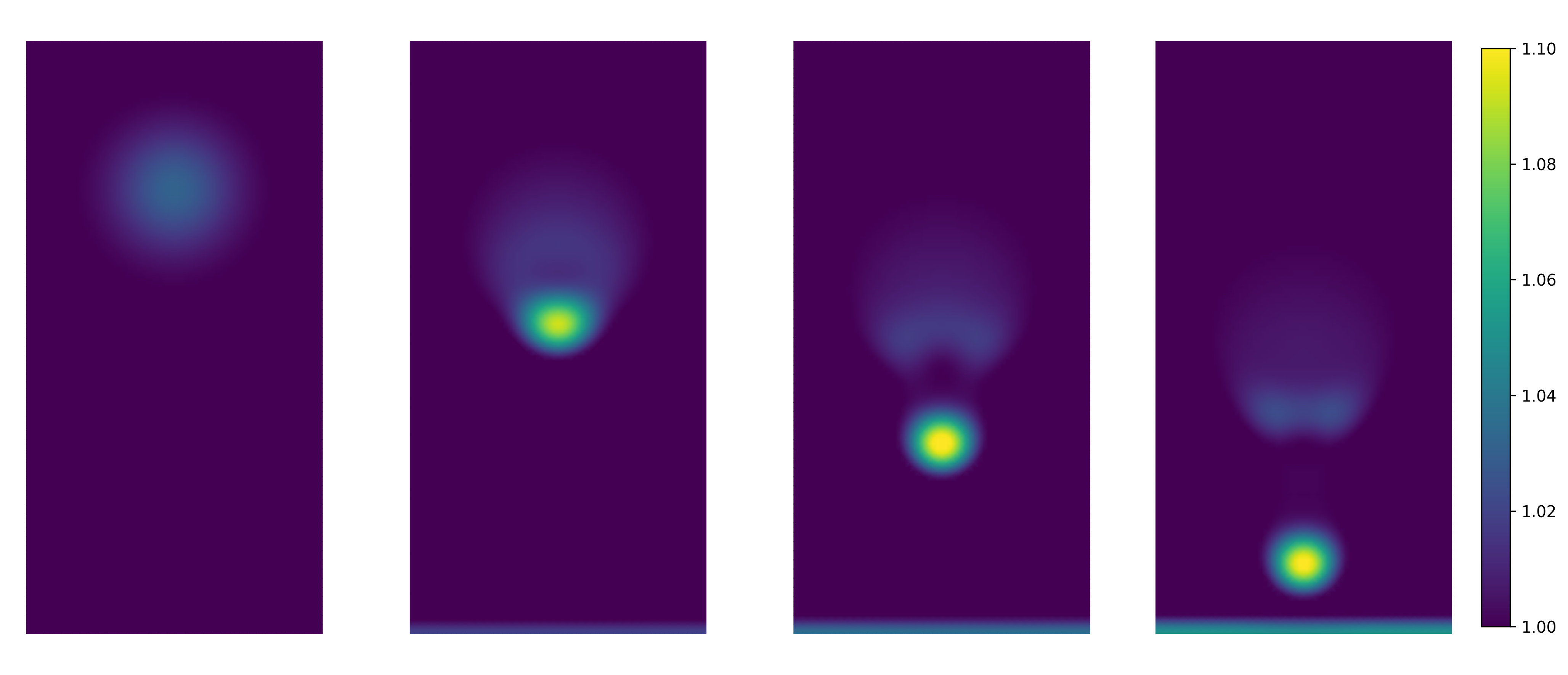}
\end{minipage}
\caption{Porosity wave moving fluid downwards and displacing the (elastic) material in the \textbf{diffusive regime} at times $t=0,4/3,8/3,4$. Concentration $\con$ and on the right the determinant $J=\det\bF$. Both for time advancing from left to right.}
\label{fig:diffusing}
\end{figure}

\begin{figure}
\centering
\hfill
\begin{minipage}{0.45\textwidth}
\includegraphics[width=0.98\textwidth,trim=0.8cm 0 0 0,clip]{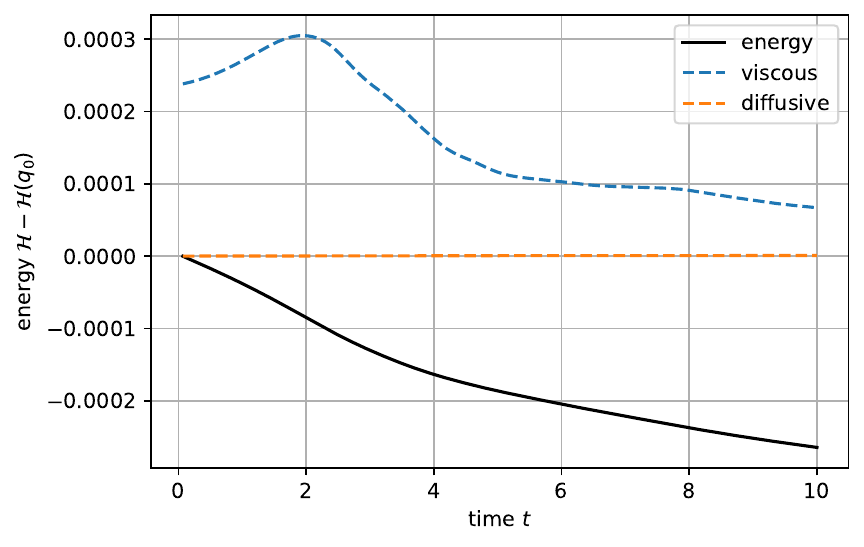}
\end{minipage}
\hfill
\begin{minipage}{0.45\textwidth}
\includegraphics[width=0.9\textwidth,trim=0.8cm 0 0 0,clip]{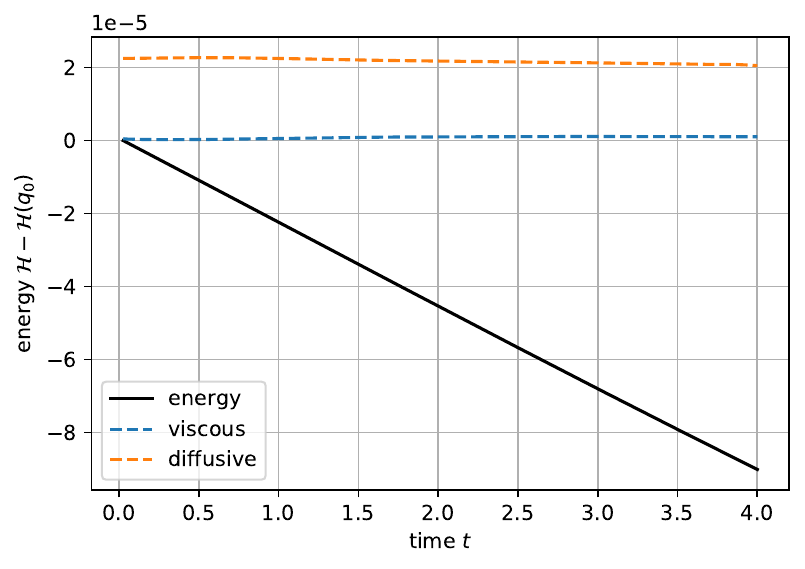}
\end{minipage}\hfill \\ 
Flowing regime. \hspace{.4\textwidth} Diffusive regime.
\caption{Decreasing system energy $\calH(\bq(t))$ (black line) and viscous dissipation $s(\partial_t \bq,\partial_t\bq)$ (blue dashed line) and diffusive energy loss $k(\eta,\eta)$ (orange dashed line) as a function of time $t$. {The energy-dissipation balance is $\tfrac{\rm d}{\mathrm{d}t}\mathcal{H}=-(s+k)$ as in \eqref{eqn:energy_descent}.}}
\label{fig:energy_flow_and_diffusion}
\end{figure}

We present in \Cref{fig:flowing} and \Cref{fig:diffusing} two examples of different types of porosity waves that are descending due to the higher mass density of the fluid, where the main difference is that the parameters give rise to a transport of fluid content primarily by convection in the former, \emph{i.e.}, the \emph{flowing regime}, and primarily by diffusion in the latter, \emph{i.e.}, the \emph{diffusive regime}. This is achieved by the different sets of parameters shown in \Cref{tab:regime_comparison}, where in the flowing regime the shear modulus, but also the diffusion constant, are both very small, whereas in the diffusive regime we suppress flow by a larger shear modulus and enhance diffusion by a larger diffusion constant. Without going into the details of the solution, this difference is most prominently observed in \Cref{fig:energy_flow_and_diffusion}, where in the flowing regime the dissipation from diffusion, encoded in $k(\bs{\eta},\bs{\eta})$, is practically negligible, whereas in the diffusive regime the Kelvin–Voigt viscous dissipation, encoded in $s(\partial_t\bs{q},\partial_t\bs{q})$, is similarly negligible. If one looks at the individual contributions (elastic, Biot, gravity, regularisation), one observes that primarily gravity $H_\textrm{grav}(c)$ changes, while the other contributions leave the porosity wave almost in local equilibrium, \emph{i.e.}, their contributions to the energy are very small.

The interpretation of these shapes is partially limited by the fact that several simplifications are applied here, \emph{i.e.}, the geometry is two-dimensional, and compaction of the solid phase due to gravity is neglected. This is both because a term $J^{-1}\varrho_0 z g_0$ is missing in the energy, and because the overall domain is fixed. The latter could be circumvented by introducing an additional indicator function for a fictitious air phase, as in \cite{peschka2025}. The main difference is that in the flowing regime, the porosity wave tends to grow in extent and leaves a more pronounced wake behind, whereas in the diffusing regime, with degenerate diffusion $D(c)=c^2 D$, the porosity wave tends to develop a more compact support. In the next example we are also going to include the compaction of the solid in a solid-fluid system, so that the resulting Biot pressure will drive porosity waves upwards.

\section{Porosity waves model}
\label{sec:poro}
 
We proceed by extending the previous single phase porous media model to a two-phase model or, more specifically, a \emph{quasi-static Euler-Euler mixture model with solid and fluid momentum balance}. 
A precise Lagrangian-Eulerian transformation as in Section \ref{sec:elasticity} is, in our opinion, not possible or practical within this setting due to the nonlocality in the Darcy dissipation. In fact, for the Lagrangian variational formulation of a two-phase system a coordinate frame has to be fixed; however every quantity which is not ''naturally`` defined in that coordinate frame will require a push-forward and pull-back mapping, for example $\flowmap_{\frm}^{-1}\circ\flowmap_{\srm}$. 

{The model that we present in this section may be regarded as a nonlinear extension of \cite{mckenzie1984generation}, incorporating hyperelastic constitutive behaviour of the solid matrix and diffusive mechanisms. Related extensions have appeared in the literature, typically emphasising more elaborate rheological descriptions, for example of Maxwell type \cite{Yarushina2015} or in combination with plastic effects and shear stresses \cite{yarushina2020model}, or by introducing shock waves in the limit of zero viscosity \cite{spiegelman1993flow}. Because of the novelties introduced in our model, {\emph{i.e.}, nonlinear elasticity based on an energetic formulation}, reliable estimates of the associated parameters are difficult to identify. We therefore adopt parameter values consistent with those commonly employed in previous studies and focus on qualitative behaviour rather than precise quantitative calibration; see, for instance, \cite{scott1984magma,wiggins1995magma} for representative length and time scales.}

When dealing with mixture models, the kinematics description is an essential ingredient. In this work, we tailor the presentation and the model specifically and assume that to each phase $\irm\in\lbrace\srm,\frm\rbrace$ there is an associated flow map $\flowmap_{\irm}$. This strategy is commonly used in the literature and for more details about specific constructions and comparisons between different approaches to multiphase kinematics and mixture flows we refer to, \emph{e.g.},   \cite{de2012theory,ehlers2009challenges,coussy2004poromechanics}.

We consider a system formed both by a solid ($\srm$) and a fluid phase ($\frm$), whose amount is locally given by a volume fraction $\phi_{\irm}:\Omega\to\R$ for $\irm \in \{\srm,\frm\}$. This amount is normalised to sum up to 1, i.e., $\phi_\frm(\bx) + \phi_\srm(\bx)=1$ for all $\bx\in\Omega$, which allows us to express the solid volume fraction in terms of the fluid one. The state vector $\bq$ is composed by the Eulerian displacement of the solid and the fluid part, cf. Section \ref{sec:kinematics}, and by the fluid content $\zeta$, i.e.,
\begin{align*}
    \bq = \begin{pmatrix}
        \bu_{\srm}:\Omega\to\R^d \\ \bu_{\frm}:\Omega\to\R^d \\ \zeta:\Omega\to\R
    \end{pmatrix}\in \calQ\,.
\end{align*}
The fluid content $\zeta$ was introduced by Biot \cite{biot1941} in the description of a porous elastic material and it relates to the fluid mass density $\rhof:\Omega\to\R$ and the fluid volume fraction or \emph{porosity} $\phi_{\frm}:\Omega\to\R$ through
\begin{align}\label{eqn:def-fluid-content}
    \zeta(t,\bx):= \frac{\phi_{\frm}(\bx) \rhof(\bx)}{\Jf(t,\bx)}\,,
\end{align}
where we recall that $\bFalei=(\idty - \nabla\bu_{\irm}(t,\bx))^{-1}$, $\Ji=\det\bFalei$ and the velocity field $\bvi = \bF_{\irm}\partial_t \bu_{\irm}$ for $\irm \in \lbrace\srm,\frm\rbrace$. We highlight that \eqref{eqn:def-fluid-content}, and the corresponding definitions of $\bu_{\srm}, \bu_{\frm}$, define a transformation $\trafo:\bcalQ\to\calQ$ in the spirit of \eqref{eqn:def-trafo-el-di}. Upon differentiation of \eqref{eqn:def-fluid-content}, we get the following transport equation for $\zeta$
    \begin{align*}
        \partial_t \zeta(t,\bx) + {\rm div}\left(\zeta(t,\bx)\bv_{\frm}\right) = 0\,.
    \end{align*}

Additionally, we emphasise that this two-phase model is characterised by two distinct flow maps or displacements. This choice allows us to precisely transport each individual material component, which comes at the price of higher model complexity and requires additional material parameters to be inferred from experimental data.
\par
\subsection*{Weak and strong formulation}
We define the three bilinear forms characterising \eqref{eqn:MAUS}: Consider a state $\bq\in\calQ$, $\bv, \bv_1, \bv_2\in\calV$ and $\bw,\bw_1,\bw_2\in\calW$ and set
\begin{subequations} \label{eqn:porosity-bilinear-forms}
    \begin{align}
        \brm(\bq;\bv,\bw)&:=\int_{\Omega}    \bw_{\zeta}\bv_{\zeta} + \bw_{\zeta}\nabla\cdot(\zeta\bF_{\!\frm}\bv_{\bu_\frm})\,{\rm d}\bx\,,
        \\ \label{eqn:diff-diss}
        \krm(\bw_1,\bw_2)&:= \int_{\Omega}\nabla\bw_{1,\zeta}\cdot D(\zeta)\nabla\bw_{2,\zeta}\,{\rm d}\bx\,,
        \\
        \srm(\bq;\bv_1,\bv_2)&:= \srm_{\rm Darcy}(\bq;\bv_1,\bv_2) + \srm_{\rm visc}(\bq;\bv_1,\bv_2)\,,
        \\\label{eqn:darcy-diss}
        \srm_{\rm Darcy}(\bq;\bv_1,\bv_2) &:= \int_{\Omega}M(\zeta)\left(\bF_{\!\srm}\bv_{1,\bu_\srm} - \bF_{\!\frm}\bv_{1,\bu_\frm}\right)\cdot\left(\bF_{\!\srm}\bv_{2,\bu_\srm} - \bF_{\!\frm}\bv_{2,\bu_\frm}\right)\,{\rm d}\bx\,,
        \\\nonumber
        &\text{with } M(\zeta):=M_0\left(\frac{\zeta_0}{\zeta}\right)^2\,\, \text{ and }\,\, M_0>0\,,
        \\ \label{eqn:visc-diss}
        \srm_{\rm visc}(\bq;\bv_1,\bv_2) &:= \int_{\Omega} \nu_{\srm}\nabla\left(\bF_{\!\srm}\bv_{1,\bu_\srm}\right){:}\nabla\left(\bF_{\!\srm}\bv_{2,\bu_\srm}\right) + \nu_{\frm}\nabla\left(\bF_{\!\frm}\bv_{1,\bu_\frm}\right){:}\nabla\left(\bF_{\!\frm}\bv_{2,\bu_\frm}\right)\,{\rm d}\bx\,.
    \end{align}
    The bilinear form $\brm$ is the analogue of \eqref{eqn:eul-b} for the fluid content $\zeta$ being transported by the fluid velocity $\bF_\frm\bv_{\bu_\frm}$, as suggested by its definition \eqref{eqn:def-fluid-content}.
    The dissipative bilinear form $\krm$ in \eqref{eqn:diff-diss} defines a $H^{-1}$-type dissipation as in \eqref{eqn:diff_kblf} with $D(\zeta)=D_0 \zeta$ for some $D_0>0$, while \eqref{eqn:visc-diss} introduces viscous effects of Stokes-type and it is similar to the Stokes dissipation from Remark \ref{rem:overdamped-ham} with $\nu_\srm,\nu_\frm>0$ {shear} viscosity coefficients. {We mention that extensions to the operator $\srm$ to include bulk viscosity effect can be made by adding volumetric terms of the type $\nabla\cdot(\bF_{\!\srm}\bv_{i,\bu_\srm})$, see e.g. \cite{zafferi2022generic}.}
    The Darcy component \eqref{eqn:darcy-sblf} exhibits a degeneracy similar to the diffusive component $k$ that will produce solutions with compact support, cf. \eqref{eqn:diff_kblf}. In this case,  the coefficient behaves differently for $\zeta\to 0^+$ (\emph{zero-porosity limit}) since $M(\zeta)\to\infty$ while in the diffusive case we have $D(c)\to 0$ as $c\to 0^+$. Note, the form $s_\textrm{Darcy}$ would be nonlocal in Lagrangian coordinates.
\end{subequations}

\begin{subequations}\label{eqn:porosity_energy}
    The free energy of the system is the sum of four different contributions: an elastic solid, a gravitational, a coupling or Biot term and a regularising term
    \begin{align}
        \calH(\bq):=&\int_{\Omega}H_{\rm elastic}(\bF_{\!\srm}) - H_{\rm gravity}(\Js,\zeta) + H_{\rm Biot}(\Js,\zeta)+ H_{\rm reg}(\zeta,\nabla\zeta)\,\dbx\,,
        \\  
        H_{\rm elastic}(\bF_{\!\srm}) :=& \frac{\mu}{2}{\rm tr}\left(\bC_{\rm iso}-\idty\right)
        \\ \label{eqn:gr_energy}
        H_{\rm gravity}(\Js,\zeta) :=&\, gz \left(\zeta + \frac{\rho_s (1-\phi)}{\Js}\right)\,,
        \\ \label{eqn:mix_energy}
        H_{\rm Biot}(\Js,\zeta) :=& \frac{\kappa}{2}\pre^2\,, \qquad \pre:=\left((\zeta - \zeta_0) - (\Js - 1)\right)\,,
        \\ \label{eqn:reg_energy}
        H_{\rm reg}(\zeta,\nabla\zeta) =& \frac{\epsilon_1}{2}\vert\nabla\zeta\vert^2 + \epsilon_2\zeta\left(\log\left(\zeta/{\zeta^*}\right)-1\right)\,,
    \end{align}
    where $\mu,\kappa>0$ are material parameters, $\pre$ is the pore pressure, $\epsilon_1,\epsilon_2>0$ are weights for the regularisation effects and, as before, with coordinate vector {$\bx=(x,z)\in\R^d$ where $x\in\R^{d-1}$ and $z\in\R$, $\bs{g} = (0,g)$ where $0\in\R^{d-1}$ and $g\in\R$}.
    The nonlinear elastic energy density is unchanged from \eqref{eqn:energy_single_phase}, while for the pressure and gravitational energy densities adjustments were made for the quasi-static Euler-Euler mixture model. In particular, the Biot term couples now changes in the solid volume $\Js$ with the fluid content $\zeta$, where pressure defined in \eqref{eqn:mix_energy} is actually a special case of the usual pore pressure for poroelastic materials, cf. \cite{mielke2013homogenization}, 
    $\pre = \zeta-\mathbb{A}:(\bF_{\!\srm}- \idty)$ defined by introducing the tensor of the so-called \emph{Biot's coefficients}, where we have set $\mathbb{A}$ to be the determinant. This modelling choice produces a volumetric-isochoric split in the energy contributions, where only solid volumetric changes influence the fluid content and conversely, but shear deformations are still allowed by the elastic energy. As we shall see in Section \ref{sec:sim_pw}, this strategy will allow us to capture porosity waves and large deformations while keeping the model relatively simple. Gravity has been extended to include both phases and, by comparing \eqref{eqn:gr_energy} with \eqref{eqn:def-fluid-content}, one can see that both mass densities are accounted for; in fact, by tuning the ratio between $\rhof$ and $\rho_s$, one is able to create ascending or descending porosity waves.  Variations of these two energy densities will be the main contributors to the driving force of the system. 
    The last energy density encodes some higher regularising term in $\zeta$, as in the last term in \eqref{eqn:hamiltonian-eul-conc}, and a logarithmic contribution that prevents unphysical solutions, \emph{i.e.} $\zeta<0$.
    \end{subequations}
    From \eqref{eqn:porosity_energy} we can compute the thermodynamic driving forces for any $\bv=(\bv_{\bu_\srm},\bv_{\bu_\frm},\bv_{\zeta})\in\calV$
    \begin{align}\label{eqn:pw-driving-force}
    \begin{split}
         \left\langle\D\calH(\bq),\bv\right\rangle = \int_{\Omega} \left(\bF_{\!\srm}^{\top}\partial_{\bF_{\!\srm}}(H_{\rm elastic}(\bF_\srm))\bF^{\top}_{\srm} + \kappa \pre \Js \bF_{\!\srm}^{\top} - gz\frac{\bF_{\!\srm}^{\top}}{\Js} \right):\nabla\bv_{\bu_{\srm}}\dbx
        \\  + \int_{\Omega}\left(\kappa \pre + \epsilon_2\log(\zeta/\zeta^*)\right)\bv_{\zeta} + \epsilon_1 \nabla\zeta\cdot\nabla\bv_{\zeta} + gz\bv_{\zeta}\,\dbx = {(\bs{f}_{\srm},\bv_{\bu_\srm}) + (\bs{f}_{\zeta},\bv_{\zeta})}\,,
    \end{split}
    \end{align}
{where we recall that $(\cdot,\cdot)$ denotes the $L^2$ inner product} and that variations of $\bu_f$ do not produce any driving force, thus $\bs{f}_{\frm}\equiv0$. This is to be expected since only variations in $\Jf$ can generate thermodynamic potentials, but this is included in the definition of fluid content \eqref{eqn:def-fluid-content} thus not appearing explicitly in the energy. 
The disappearance of the determinant $\Jf$ is a common reduction strategy for fluid systems, see \emph{e.g.} \cite{morrison1998hamiltonian,zafferi2022generic}. 
The functionals defined in \eqref{eqn:pw-driving-force} and the bilinear forms \eqref{eqn:porosity-bilinear-forms} fully determine a weak formulation of a PDE system through \eqref{eqn:MAUS} as follows: 
\begin{subequations}\label{eqn:pw-weak-system}
    \begin{align}
        \label{eqn:pw-weak-system-a}\int_{\Omega}\bw_{\zeta}\left(\partial_{t}\zeta + \nabla\cdot(\zeta\bF_{\!\frm}\partial_t\bu_{\frm})\right)\dbx =& -\int_{\Omega}\nabla\bw_{\zeta}\cdot D(\zeta)\nabla\eta\,{\rm d}\bx\,,
        \\ \label{eqn:pw-weak-system-b}
        \int_{\Omega}\nu_{\srm}\nabla\left(\bF_{\!\srm}\partial_t\bu_\srm\right){:}\nabla\left(\bF_{\!\srm}\bv_{\bu_\srm}\right) &-M(\zeta)\left(\bF_{\!\srm}\dot{\bu}_{\srm} - \bF_{\!\frm}\dot{\bu}_{\frm}\right)\cdot\bF_{\!\srm}\bv_{\srm} \,\dbx
        \\ \nonumber =&\int_{\Omega}\left(\bF_{\!\srm}^{\top}\partial_{\bF_{\!\srm}}(H_{\rm elastic}(\bF_\srm))\bF^{\top}_{\srm} + \kappa \pre \Js \bF_{\!\srm}^{\top} -gz \frac{\bF_{\!\srm}^{\top}}{\Js}\right):\nabla\bv_{\bu_{\srm}}\dbx\,,
        \\\nonumber
        \int_{\Omega}\nu_{\frm}\nabla\left(\bF_{\!\frm}\partial_t\bu_{\frm}\right){:}\nabla\left(\bF_{\!\frm}\bv_{\bu_\frm}\right) &+ M(\zeta)\left(\bF_{\!\srm}\dot{\bu}_{\srm} - \bF_{\!\frm}\dot{\bu}_{\frm}\right)\cdot\bF_{\!\srm}\bv_{\frm} \,\dbx
        \\ \label{eqn:pw-weak-system-c}
        =&\int_{\Omega} \eta \nabla\cdot(\zeta\bF_{\!\frm}\bv_{\frm})\dbx\,,
        \\ \label{eqn:pw-weak-system-d} 
        \int_{\Omega}\eta\bv_{\zeta}\dbx =& \int_{\Omega}\left(\kappa \pre + \epsilon_2\log(\zeta/\zeta^*) + gz\right)\bv_{\zeta} + \epsilon_1 \nabla\zeta\cdot\nabla\bv_{\zeta}\,\dbx\,.
    \end{align}
Equation \eqref{eqn:pw-weak-system-a} conservation of fluid content, while \eqref{eqn:pw-weak-system-d} sets the chemical potential. The last two equations \eqref{eqn:pw-weak-system-b} and \eqref{eqn:pw-weak-system-c} are the solid and fluid momentum balance without inertia contributions.
\end{subequations}
The strong form can be derived by collecting and integrating by parts the terms corresponding to:
\begin{subequations}\label{eqn:pw-strong-system}
    \begin{align}
        \bw_{\zeta}&: \dot{\zeta} + \nabla\cdot\left(\zeta\bF_{\!\frm}\dot{\bu}_{\frm}\right)&&=\nabla\cdot\left(D(\zeta)\nabla\bs{\eta}\right)\,,
        \\
        \bv_{\zeta}&: \bs{\eta} &&= \bs{f}_{\zeta}\,,
        \\
        \bv_{\bu_\srm}&: -M(\zeta)\bF_{\!\srm}^{\top}\left(\bF_{\!\srm}\dot{\bu}_{\srm} - \bF_{\!\frm}\dot{\bu}_{\frm}\right) +   \bF_{\!\srm}^{\top}\nabla\cdot\left(\nu_{\srm}\nabla(\bF_{\!\srm}\dot{\bu}_{\srm})\right)&&=\bs{f}_{\srm}\,,
        \\
        \bv_{\bu_\frm}&: M(\zeta)\bF_{\!\frm}^{\top}\left(\bF_{\!\srm}\dot{\bu}_{\srm} - \bF_{\!\frm}\dot{\bu}_{\frm}\right) + \bF_{\!\frm}^{\top}\nabla\cdot\left(\nu_{\frm}\nabla(\bF_{\!\frm}\dot{\bu}_{\frm})\right) - \bF_{\!\frm}^{\top}\nabla \bs{\eta}&&=0\,.
    \end{align}
\end{subequations}

\subsection{Discretisation}
Following the same discretisation recipe introduced before, we are going to discretise \eqref{eqn:pw-weak-system} on the level of the involved bilinear forms and apply a semi-implicit scheme in time and finite elements in space, \emph{i.e.}, the derivative of the free energy $\D\calH(\bq)$ is discretised implicitly while the state dependence in the bilinear forms $\brm, \krm, \srm$ is discretised explicitly in $\bq$. We assume that we have a sequence of time-steps $0=t^0<t^1<\ldots<t^N=T$ and write for every $k\in\lbrace 0,\ldots,N\rbrace$ $\bq^k(\bx)=\bq(t^k,\bx)$ as well as $\bs{\eta}^k(\bx)=\bs{\eta}(t^k,\bx)$ or also simply $\bq^k=\bq(t^k)$ and $\bs{\eta}^k=\bs{\eta}(t^k)$. We approximate the time-derivative by $\partial_t\bq^k=(\bq^k - \bq^{k-1})/\tau^k$ with $\tau^k=t^{k}-t^{k-1}$. Consider a given triangulation $\mathcal{T}_h$ of the domain $\Omega = \cup_{T\in\mathcal{T}_h}T$ and the generic finite-dimensional function spaces defined in \eqref{eqn:FEM-space}. Then, we define our approximate states via, cf. \eqref{eqn:approx-state}:
\begin{align*}
    &\bq^k = (\bu^k_{\srm},\bu^k_{\frm},\zeta^k)\in V_h^{d,2}\times V_h^{d,2}\times V_h^{1,1}=:\calV_h\,, &&\bs{\eta}^k\in V_h^{1,1}=:\calW_h\,.
\end{align*}
We complement the spaces for $\bu^k_{\srm}, \bu^k_{\frm}$ with Dirichlet boundary conditions $\bu^k_{\srm}\cdot\bs{\nu}=\bu^k_{\frm}\cdot\bs{\nu}=0$ on the sides and $\bu^k_{\srm}=\bu^k_{\frm}=0$ on top and bottom. For the fluid content $\zeta^k$ we adopt homogeneous Neumann boundary conditions, as well as for the chemical potential $\bs{\eta}^k$. Finally, for a given $\bq^{k-1}$ we solve the nonlinear finite dimensional system \eqref{eqn:DMAUS} for $\bq^k$ and $\bs{\eta}^k$. 

\subsection{Simulations}\label{sec:sim_pw}
\begin{figure}[t]
\centering
\begin{minipage}{.49\textwidth}
\includegraphics[width=\textwidth, trim=4.5cm 0 2cm 0,clip]
{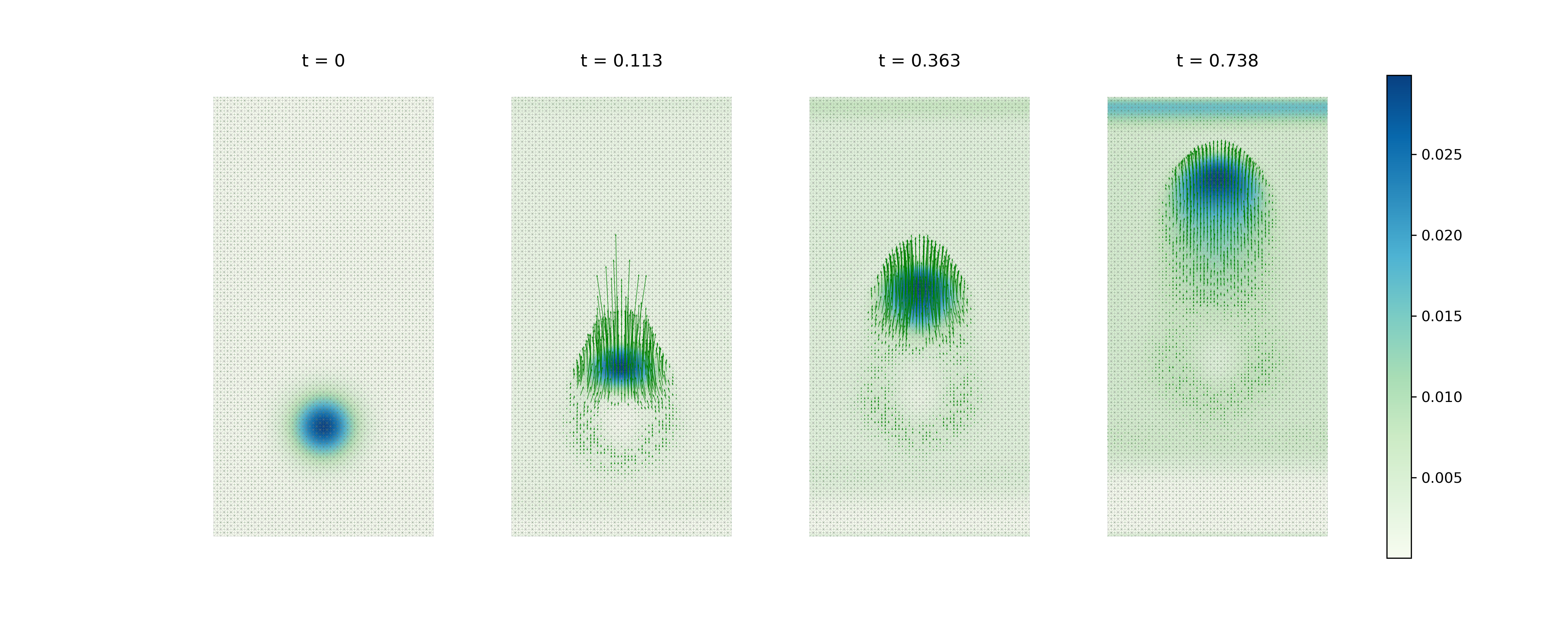}%
\end{minipage}
\begin{minipage}{.49\textwidth}
\includegraphics[width=\textwidth ,trim=4.5cm 0 2cm 0,clip]
{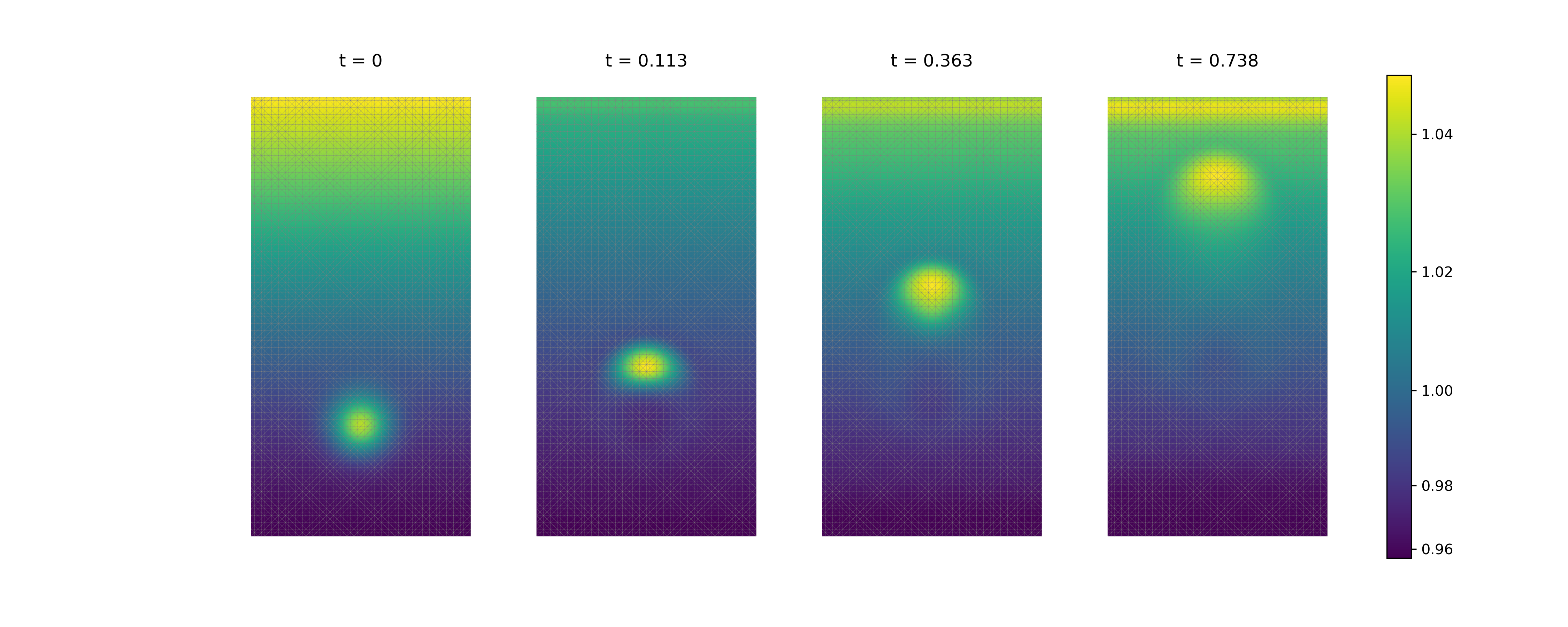}
\end{minipage}
\\
\caption{Porosity wave moving fluid upwards and displacing the (elastic) material. On the left the rescaled fluid content $\zeta$ with green arrows denoting fluid velocity are illustrated. On the right the determinant $J_s=\det\bF_s$. Both for time advancing from left to right.}
\label{fig:compaction_poro}
\end{figure}

\begin{figure}[t]
    \centering
    \includegraphics[width=0.75\linewidth]{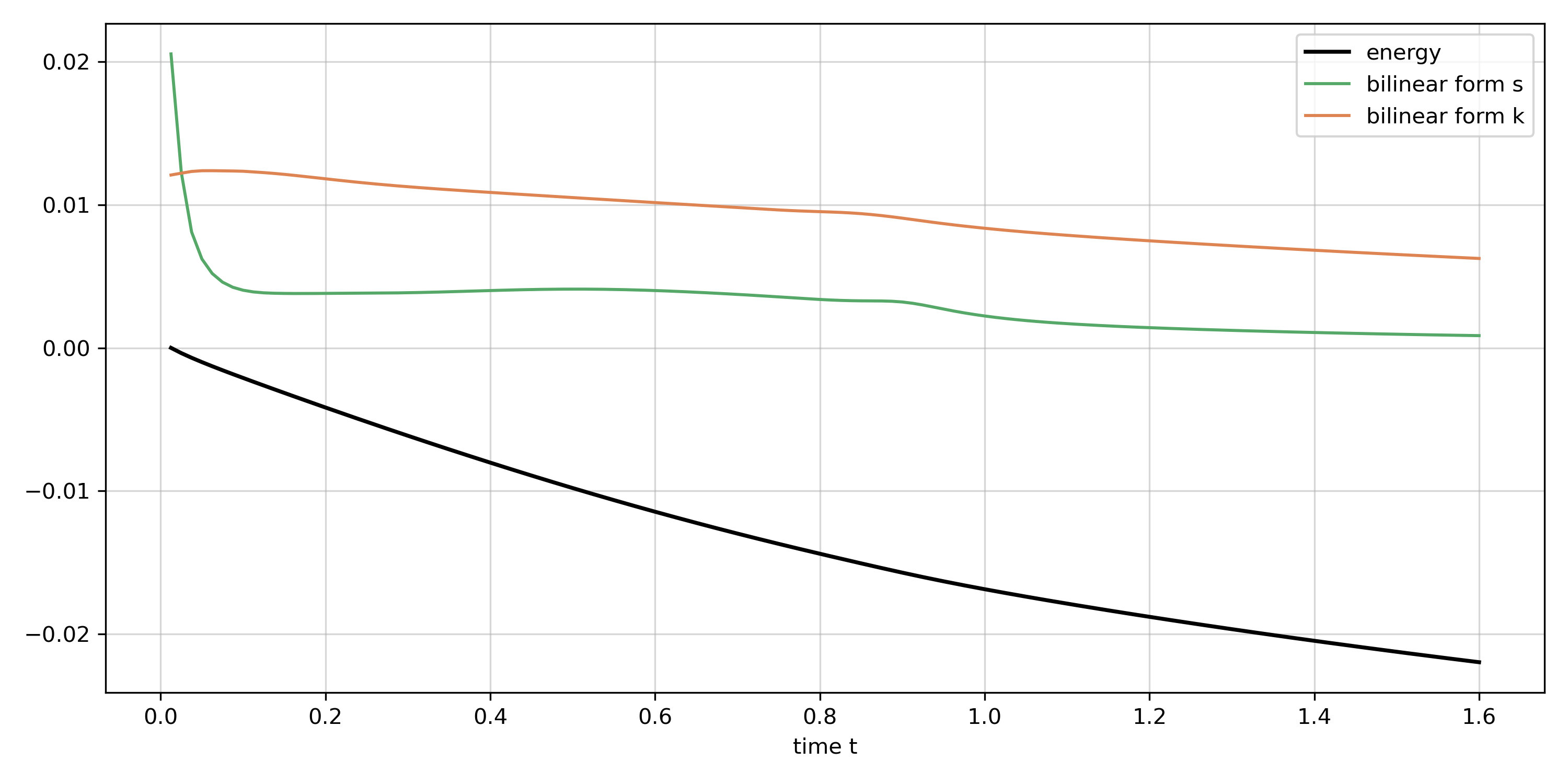}
    \caption{Relative system Energy $\calH(\bq(t))-\calH(\bq(0))$ (black line), bilinear form $s(\partial_t \bq,\partial_t \bq)$ (blue) and $k(\bs{\eta}(t),\bs{\eta}(t))$ behaviour (orange) as function of time. The energy-dissipation balance is $\tfrac{\rm d}{\mathrm{d}t}\mathcal{H}=-(s+k)$ as in \eqref{eqn:energy_descent}.}
    \label{fig:total_energy}
\end{figure}

For the simulations, we consider a 2-dimensional rectangular domain $\Omega = (0,L)\times(0,2L)$ with components $\bx=(x,z)\in\Omega$. In order to simulate rising porosity waves, we start from a localised excess of fluid content expressed by a Gaussian distribution as in \eqref{eqn:inital_datum1} 
\begin{align}
    \zeta^0(\bx) = \bar{\zeta} + c\exp(-\tfrac{1}{r^2}\vert \bx - \bx_0\vert^2)\,,
\end{align}
where we choose $\bx_0 = (\tfrac{L}{2},\tfrac{L}{2})$ and $\bar{\zeta}$ is usually referred to as background fluid content or porosity. 
The initial data for the displacements $\bu_\srm^0$ and $\bu_\frm^0$ is then obtained as stationary points to the Hamiltonian $\calH=\calH(\bq)$ defined in \eqref{eqn:porosity_energy} with $\zeta=\zeta^0$. 
We normalize the equations \eqref{eqn:pw-strong-system} by rescaling the domain with the length $L = 10^4\,m$, the hydrostatic pressure $P_0=\rhof g L$ with $g = 10\,m/s^2$ (an approximation of the gravitational constant) and $\rhof = 10^3\, kg/m^3$ (an approximation of pure H$_2$O fluid) and a typical time length $T=10^{11}\,s$, see, e.g., the examples in \cite{connolly2012analytical} for more details about the other parameters. Lastly, we have used the low porosity approximation, \emph{i.e.}, $1-\phi \approx 1$ to simplify the expression for the gravitational energy. This allows us to avoid introducing a $1/\Js$ term in the gravitational energy (since $\rho_s(1-\phi)$ is transported by the solid velocity) which might generate stability issues within the code. After renormalisation, the variable $\zeta$ lies within the values of the porosity $\phi$, cf. \eqref{eqn:def-fluid-content}, thus we set $\bar{\zeta}=5\times10^{-3}$ and $\max_{\Omega}{\zeta^0}=5\times10^{-2}$. In Figure \ref{fig:compaction_poro}, we present time snapshots of the numerical solutions to \eqref{eqn:pw-weak-system}, illustrating the characteristic compaction–decompaction dynamics of the solid matrix, whose motion drives the transport of excess fluid. In fact, one can observe the solid matrix expanding $(\Js>1)$ at the top and around the fluid blob and compressing $(\Js<1)$ at the bottom. 
In contrast with the simulations of Section \ref{sec:pw-single-phase}, here we are capable of capturing this mechanism thanks to the employment of a flow map for each phase and the subsequent gravitational energy density \eqref{eqn:gr_energy}. 
In this case, the ratio in the mass densities $\rho_s/\rhof$ pushes the fluid to rise and the degeneracy in the Darcy coefficient (together with the low background fluid content $\bar{\zeta}$) forces the fluid to localize and to travel compactly instead of diffusing. Lastly, the pore pressure arising from the Biot term forces the solid to expand whenever there is an excess fluid content, \emph{i.e.}, $\zeta > \bar{\zeta}$. The ratio $\rho_s/\rhof$ strongly influences the propagation speed of the porosity wave. 
To show that the thermodynamic consistency of the model induced by our variational framework, cf. \eqref{eqn:energy_descent}, is preserved after discretisation we plot the energy and the dissipative bilinear forms in Figure \ref{fig:total_energy}.
The use of the Stokes-Stokes and the diffusion dissipations mechanisms completes the description delivering effectively a porosity and not a shock wave, see, \emph{e.g.}, \cite{spiegelman1993flow}.  

\begin{table}[htbp]
\centering
\caption{Values of rescaled parameters of the simulation for Figure \ref{fig:compaction_poro}. The system was non-dimensionalised with the hydrostatic pressure $P_0=\rhof g L=10^8\,Pa$ and time length of $T=10^{11}s\approx 3\times10^3\, years$} \small
\begin{tabular}{|l|c|c|c|c|c|c|c|c|c|c|}
\hline
\textbf{Parameters}& $\kappa/P_0$ & $\rho_s/P_0$ & $\epsilon_1/P_0L^2$ & $\epsilon_2/P_0$  & $\mu/P_0$  & $D_0P_0T/L^2$ & $M_0P_0T/\nu_fL^2$ & $\nu_f/P_0T$ & $\nu_s/P_0T$ 
 \\ 
 \textbf{Values} & $100$ & $4$ & $10^{-2}$ & $10^{-2}$ & $1$ & $10^{-1}$ & $1$ & $10^{-4}$ & $10^{-4}$ \\
\hline
\end{tabular}
\label{tab:parameters}
\end{table}

\section*{Conclusions and outlook}
{In this paper, Lagrangian and Eulerian descriptions of nonlinear elasticity and poroelasticity have been revisited and formulated within a unified variational framework that employs saddle-point structures and formal transformations between coordinate frames \cite{zafferi2022generic,peschka2022variational,schmeller2023gradient}. Within this framework, the paper introduces a versatile energy-based modelling approach for viscous poroelastic materials with diffusive transport in both Lagrangian and Eulerian frames. Key results are the development of refined Eulerian formulations using the reference map concept, a novel structure-preserving discretisation strategy based on saddle-point problems, and the demonstration of spatial and temporal numerical convergence. One key finding are certain closure relations that bilinear forms of mixed formulations need to satisfy in order to ensure consistency of Lagrangian and Eulerian formulations. Additionally, the paper demonstrates the existence of poroelastic waves within this class of models.

We focused on geological applications related to porosity waves, introducing and motivating step by step new conservative and dissipative effects within the variational structure. In parallel, we have leveraged these saddle-point structures to provide practical finite-element based discretisations in order to study convergence of solutions numerically.}
This family of thermodynamically consistent variational frameworks can be applied to a wide range of physical phenomena. A further step in this direction would be to apply it to other classes of rheology, \emph{e.g.}, Maxwell rheology, in contrast to the Kelvin–Voigt model used in this paper. This modification of the constitutive material laws would lead to different types of porosity waves; see, \emph{e.g.}, \cite{Yarushina2015}. The transformation strategy introduced for mapping Lagrangian to Eulerian variables can be implemented as a simple change of variables or even as an extension map. For example, polyconvex elastic energies are usually written in terms of the invariants of the strain. In that case, a transformation $\trafo$ (or an extension) mapping $\bbF$ to the invariants and the inclusion of the kinetic energy to generate convective derivatives could be employed to reformulate the PDE system, potentially leading to a more structured and analytically tractable formulation, see, \emph{e.g.}, \cite{tornquist2025}.

\section*{Acknowledgement}
 AZ acknowledges the funding by the DFG-Collaborative Research Center 1114 \emph{Scaling Cascades in Complex Systems}, project \#235221301, C09 \emph{Dynamics of rock dehydration on multiple scales}. DP thanks for the funding within the DFG Priority Program SPP 2171 \emph{Dynamic Wetting of Flexible, Adaptive, and Switchable Surfaces}, project \#422792530. { Both authors acknowledge insightful discussions about gradient flows and mathematical models in geoscience with Marita Thomas and Tom{\'a}{\v{s}} Roub{\'\i}{\v{c}}ek, as well as valuable explanations and geological insights on porosity wave provided by Timm John and Johannes C. Vrijmoed.}

Generative AI tools (ChatGPT version GPT-5.1 and GitHub Copilot) were used solely for language editing (grammar, spelling) and programming support (documentation, debugging, visualization). All AI-assisted content was checked carefully by the authors. No AI tools were used to generate original research results (design, interpretation, conclusions).

\bibliographystyle{unsrt}
\bibliography{references}
\end{document}